\documentclass{article}%
\usepackage{amsmath}
\usepackage{amssymb}
\usepackage{amsfonts}
\usepackage{graphicx}%
\providecommand{\U}[1]{\protect\rule{.1in}{.1in}}
\newtheorem{theorem}{Theorem}

\newtheorem{definition}[theorem]{Definition}

\newtheorem{lemma}[theorem]{Lemma}

\newtheorem{remark}[theorem]{Remark}

\newenvironment{proof}[1][Proof]{\noindent\textbf{#1.} }{\ \rule{0.5em}{0.5em}}
\begin{document}

\title{On forward attractors for nonautonomous dynamical systems with application to
the asymptotically autonomous Chafee-Infante equation}
\author{Jos\'{e} Valero\\{\small Centro de Investigaci\'{o}n Operativa }\\{\small Universidad Miguel Hern\'{a}ndez de Elche}\\{\small Avda. de la Universidad s/n}\\{\small 03203-Elche (Alicante), Spain}\\{\small E.mail:\ jvalero@umh.es}}
\date{}
\maketitle

\begin{abstract}
In this paper, for nonautonomous dynamical systems, we give first general
conditions ensuring that a pullback attractor is a forward attractor as well
in both the single and multivalued frameworks. In particular, we consider
asymptotically autonomous systems. After that, and this is the main result of
this paper, we apply these abstract theorems to the asymptotically autonomous
Chafee-Infante equation. Finally, applications to ordinary and parabolic
differential inclusions are given.

\end{abstract}

\bigskip



\bigskip

\section{Introduction}

The theory of pullback attractors for (single or multivalued) nonautonomous
dynamical systems has been succesfully developed and, as in the autonomous
situation, general dissipative and compactness assumptions ensure their
existence along with some good properties (see e.g. \cite{CarKlo09},
\cite{CLMV03}, \cite{CarLukReal}, \cite{CLRBook}, \cite{KapKasVal2011},
\cite{Schmalfuss00}). Roughly speaking, in this theory the convergence of the
solutions of the system at a final fixed final time $t$ while the initial time
$t_{0}$ tends to $-\infty$ is considered.

A different task appears when we fix the initial time $t_{0}$ and study the
convergence of the solutions of the systems as the final time $t$ tends to
$+\infty$. In this case, two types of attractors have been defined in the
literature: the uniform atractor and the forward attractor. For uniform
attractors there exists a well developed theory ensuring their existence (see
e.g. \cite{ChepVishik}, \cite{ChepVishikBook}, \cite{MV00}, \cite{KKVZBook}).
However, the situation is quite different if we speak about forward
attractors. Roughly speaking, in the theory of uniform attractors a compact
set attracting the solutions as the final time goes to $+\infty$ is
considered, while in the case of forward attractors the solutions must
approach a family of compact sets (parametrized by time) as time goes to
$+\infty$. Unfortunately, general dissipative and compactness assumptions do
not guarantee the existence of forward attractors.

In this paper, we give sufficient conditions ensuring that a pullback
attractor is a also a forward attractor. Basically, if $\mathcal{A}%
=\{\mathcal{A}(t)\}$ is a pullback attractor for a process $U$, we define its
$\omega_{0}$-limit set by
\[
\omega_{0}\left(  \mathcal{A}\right)  =Liminf_{t\rightarrow+\infty}%
\mathcal{A}(t)=\{y\in X:\lim_{t\rightarrow+\infty}dist(y,\mathcal{A}(t))=0\}
\]
and prove that if every forward $\omega$-limit set $\omega_{f}(B,t_{0})$ of
the forward asymptotically compact process $U$, where $B$ is a bounded set of
the phase space and the initial time $t_{0}$ is arbitrary, belong to
$\omega_{0}\left(  \mathcal{A}\right)  $, then $\mathcal{A}$ is a forward
attractor. In the specific situation where the process $U$ is asymptotically
autonomous, it is shown that this condition is satisfied if the pullback
attractor is continuous with respect to the limit autonomous attractor as time
tends to $+\infty$. We prove these results for both single and multivalued processes.

In \cite{ChebKloSch}, \cite{WangLiKloeden} for a cocycle $\phi$ depending on a
parameter $p\in P$ (with $P$ compact) it is proved that under some conditions
(involving the lower semicontinuity of the pullback attractor $\mathcal{A}%
=\{\mathcal{A}_{p}\}$ with respect to the parameter $p$) the pullback
attractor $\mathcal{A}$ is a uniform forward attractor as well. This result is
applied to examples in which the nonautonomous term is periodic. In
\cite[Theorem 3.4]{KloedenLorenz} (see also \cite{CuiKloeden}) a sufficient
and necessary condition involving the $\omega$-limit set of $\mathcal{A}$ is
used to prove that the pullback attractor is a forward one. However, we have
doubts about the correctness of the proof of that result. In
\cite{CuiFigLangaNasc}, in the single-valued situation, it is proved that if
the minimal attractor of $U$, given by $\mathcal{A}_{\min}=\overline
{\underset{B\text{ bounded}}{\cup}\underset{t_{0}\in\mathbb{R}}{\cup}%
\omega_{f}(B,t_{0})}$, satisfies that $dist\left(  \mathcal{A}_{\min
},\mathcal{A}\left(  t\right)  \right)  \overset{t\rightarrow+\infty
}{\rightarrow}0$, then $\mathcal{A}$ is a forward attractor. We show that this
condition is equivalent to our assumption $\omega_{f}(B,t_{0})\subset
\omega_{0}\left(  \mathcal{A}\right)  .$

We apply the abstract results to three problems:

\begin{itemize}
\item The Chafee-Infante equation:%
\[
\left\{
\begin{array}
[c]{l}%
\dfrac{\partial u}{\partial t}-\dfrac{\partial^{2}u}{\partial x^{2}}=\lambda
u-b(t)u^{3},\text{ on }(\tau,\infty)\times(0,\pi),\\
u(t,0)=u(t,1)=0,\\
u(\tau,x)=u_{\tau}(x),\text{ }x\in\left(  0,\pi\right)  ,
\end{array}
\right.
\]
where $\lambda>0$ and $b:\mathbb{R}\rightarrow\mathbb{R}^{+}$ is an uniformly
continuous and differentiable function satisfying $0<b_{0}\leq b\left(
t\right)  \leq b_{1}.$ This is one of the most popular problems in
infinite-dimensional dynamical systems. The reason is that the dynamics inside
the global attractor has been fully described in the autonomous case
\cite{Henry85}. Also, the existence of the pullback attractor is well-known
\cite{CLRBook} and its structure has been partially described in
\cite{BroCarVal}, \cite{CLRBook}, \cite{BernalVidal}. Concerning the forward
dynamics it was proved in \cite[Theorem 13.15]{CLRBook} that there is one
special positive bounded complete trajectory that forward attracts all bounded
set $B$ such there is a function $\phi_{B}\in C^{1}([0,\pi])\cap H_{0}%
^{1}(0,\pi)$ satisfying $\phi_{B}(x)>0$, for $x\in\left(  0,1\right)  $, and
$\phi_{B}\geq y\geq\varepsilon\phi_{B}$ for some $\varepsilon\in(0,1)$ and all
$y\in B.$ Moreover, this special solution attracts forward any individual
nonengative solution as time tends to $+\infty$ \cite[Theorem 13.15]{CLRBook}
(see also \cite{BernalVidal} and \cite{LRRBSV07}).

Under the additional restrictions that $1<\lambda<4$ and $b\left(  t\right)
\overset{t\rightarrow+\infty}{\rightarrow}b\in\lbrack b_{0},b_{1}]$ we prove
the continuity of the pullback attractor with respect to the limit autonomous
one when $t\rightarrow+\infty$. Using this result we establish that the
pullback attractor is a forward attractor. As far as we know, this is the
first result of this kind for this problem. For $\lambda>4$ or a general
function $b$ the question remains open.

\item An ordinary differential inclusion:%
\[
\left\{
\begin{array}
[c]{c}%
\dfrac{du}{dt}+\lambda u\in b\left(  t\right)  H(u),\ t\geq\tau,\\
u(\tau)=u_{\tau},
\end{array}
\right.
\]
where $\lambda>0$, $b:\mathbb{R}\rightarrow\mathbb{R}^{+}$ is a continuous
functions satisfying $0<b_{0}\leq b\left(  t\right)  \leq b_{1}.$ The
structure of the pullback attractor was fully described in \cite{CLV16}. We
prove now that this attractor is also a forward attractor. Also, if $b\left(
t\right)  \overset{t\rightarrow+\infty}{\rightarrow}b\in\lbrack b_{0},b_{1}]$,
then we prove the continuity of the pullback attractor with respect to the
limit autonomous one when $t\rightarrow+\infty.$ Using this we apply the
general theory from the abstract part to establish that the pullback attractor
is a forward attractor.

\item A parabolic differential inclusion:
\[
\left\{
\begin{array}
[c]{l}%
\dfrac{\partial u}{\partial t}-\dfrac{\partial^{2}u}{\partial x^{2}}\in
b(t)H(u)+\omega(t)u,\text{ on }(\tau,\infty)\times(0,1),\\
u(t,0)=u(t,1)=0,\\
u(\tau,x)=u_{\tau}(x),\text{ }x\in\left(  0,1\right)  ,
\end{array}
\right.
\]
where $b:\mathbb{R}\rightarrow\mathbb{R}^{+},$ $\omega:\mathbb{R}%
\rightarrow\mathbb{R}^{+}$ are continuous functions such that $0<b_{0}\leq
b\left(  t\right)  \leq b_{1},\ 0\leq\omega_{0}\leq\omega\left(  t\right)
\leq\omega_{1},$and $H$ is the Heaviside function given in (\ref{Heaviside}).
The structure of the pullback attractor of this problem was partially
described in \cite{CVL20}. If we consider only non-negative solutions, then
the structure was fully described in \cite{Valero2021}. We prove now that the
pullback attractor is a forward attractor as well.\ In addition, if $b\left(
t\right)  \overset{t\rightarrow+\infty}{\rightarrow}b\in\lbrack b_{0},b_{1}]$,
$\omega\left(  t\right)  \overset{t\rightarrow+\infty}{\rightarrow}\omega
\in\lbrack\omega_{0},\omega_{1}]$, then we establish the continuity of the
pullback attractor with respect to the limit autonomous one when
$t\rightarrow+\infty.$ With this result at hand, we obtain that the pullback
attractor is a forward attractor using the previous abstract theorems.
\end{itemize}

\section{Existence of forward attractors}

In this section, we obtain sufficient conditions ensuring that a pullback
attractor is a forward attractor in both the single and multivalued cases. We
consider in particular the special situation when the process is
asymptotically autonomous.

\subsection{Single-valued processes}

Let $X$ be a complete metric space with metrix $\rho$ and let $\mathbb{R}%
_{\geq}^{2}=\{\left(  t,s\right)  \in\mathbb{R}^{2}:t\geq s\}$. For
$A,B\subset X$ the Hausdorff semidistance from $A$ to $B$ is given by
\[
dist(A,B)=\sup_{y\in A}\inf_{x\in B}\rho\left(  y,x\right)  ,
\]
while the Hausdorff distance is define as%
\[
dist_{H}(A,B)=\max\{dist(A,B),dist(B,A)\}.
\]

We consider a continuous process $U:\mathbb{R}_{\geq}^{2}\times X\rightarrow
X$, which satisfies:

\begin{enumerate}
\item $U(s,s,x)=x$ for all $x\in X,\ s\in\mathbb{R};$

\item $U(t,s,x)=U(t,\tau,U(\tau,s,x))$ for all $x\in X,\ t,s,\tau\in
\mathbb{R}$, $s\leq\tau\leq t;$

\item $\left(  t,s,x\right)  \mapsto U(t,s,x)$ is a continuous map.
\end{enumerate}

A family of sets $\{B(t)\}$ is said to be bounded (closed, compact) if each
$B(t)$ is bounded (closed, compact) in $X$. It is forward asymptotically
compact if any sequence $\{y_{n}\}$, where $y_{n}\in B(t_{n})$, $t_{n}%
\rightarrow+\infty$, is relatively compact. A family $\{K(t)\}$ is invariant
if $K(t)=U(t,s,K(s))$ for any $\left(  t,s\right)  \in\mathbb{R}_{\geq}^{2}$.

\begin{definition}
A family of non-empty sets $\mathcal{A}=\{\mathcal{A}(t)\}$ is called pullback
attracting if it attracts any bounded set $B\subset X$ in the pullback sense,
which means that%
\begin{equation}
dist\left(  U(t,t_{0},B\right)  ,\mathcal{A}(t))\rightarrow0\text{, as }%
t_{0}\rightarrow-\infty\text{, for any }t\text{ fixed.} \label{PullbackAttr}%
\end{equation}
It is called forward attracting if it attracts any bounded set $B\subset X$ in
the forward sense, which means that%
\begin{equation}
dist\left(  U(t,t_{0},B\right)  ,\mathcal{A}(t))\rightarrow0\text{, as
}t\rightarrow+\infty\text{, for any }t_{0}\text{ fixed.} \label{ForwardAttr}%
\end{equation}

\end{definition}

\begin{definition}
A family of non-empty compact invariant sets $\mathcal{A}=\{\mathcal{A}(t)\}$
is called a pullback attractor if it is the minimal pullback attracting family
of closed sets. The minimality property means that if $\mathcal{C}%
=\{\mathcal{C}(t)\}$ is a closed pullback attracting family, then
$\mathcal{A}(t)\subset\mathcal{C}(t)$ for any $t\in\mathbb{R}.$
\end{definition}

\begin{definition}
A family of non-empty compact invariant sets $\mathcal{A}=\{\mathcal{A}(t)\}$
is called a forward attractor if it is forward attracting.
\end{definition}

We observe that in the last definition we have dropped the minimality
property. The reason is that in general different several forward attractors
could exist, as can be seen in the simple example of the differential equation%
\[
y^{\prime}+y=t.
\]
The solution of the Cauchy problem are given by%
\[
y\left(  t\right)  =t-1+e^{s-t}\left(  y\left(  s\right)  +1-s\right)  .
\]
The invariant family $\mathcal{A}(t)=\{t-1\}$ is both the pullback attractor
(which is minimal and unique) and a forward attractor. However, the invariant
family $\mathcal{A}_{r}(t)=\{t-1+re^{-t}\}$ is also a forward attractor for
any $r\in\mathbb{R}$. Thus, there is an infinity number of different forward
attractors and neither of them is contained in the others.

Nevertheless, it is shown in \cite[Lemma 6.2]{CuiFigLangaNasc} that forward
attractors are asymptotically equivalent, which means that if $\mathcal{A}%
_{1},\mathcal{A}_{2}$ are two forward attractors, then
\[
dist_{H}(\mathcal{A}_{1}\left(  t\right)  ,\mathcal{A}_{2}\left(  t\right)
)\rightarrow0\text{ as }t\rightarrow+\infty.
\]
We can easily see that this is true in our previous example.

In this paper, we will focus on studying when a given pullback attractor is a
forward attractor.

The convergences (\ref{PullbackAttr}), (\ref{ForwardAttr}) are equivalent for
an invariant family $\mathcal{A}$ when the convergences are uniform
\cite[Lemma 1.13]{CLRBook}, but not in general. Our aim is to establish some
sufficient conditions implying that a pullback attractor is a forward one. It
is important to point out that, while for pullback attractors there is a
general theory of their existence (see, for instance, \cite{CLRBook}), it is
much more difficult to establish that a forward attractor exists.

A complete trajectory of a process $U$ is a function $\xi:\mathbb{R}%
\rightarrow X$ such that $\xi\left(  t\right)  =U(t,s,\xi\left(  s\right)  )$
for all $t\geq s.$ A complete trajectory $\xi$ is bounded if $\cup
_{t\in\mathbb{R}}\xi\left(  t\right)  $ is a bounded set. If a pullback
attractor $\mathcal{A}$ is globally bounded, that is, $\cup_{t\in\mathbb{R}%
}\mathcal{A}(t)$ is bounded, then it can be characterized as the union of of
bounded complete trajectories \cite[Corollary 1.18]{CLRBook}, which means that%
\[
\mathcal{A}(t)=\{\xi\left(  t\right)  :\xi\text{ is a bounded complete
trajectory}\}.
\]

The forward $\omega$-limit set of a set $B\subset X$ for $t_{0}\in\mathbb{R}$
is defined by
\begin{equation}
\omega_{f}(B,t_{0})=\{y:y=\lim_{n\rightarrow\infty}y_{n}\text{, }y_{n}\in
U(t_{n},t_{0},B),\ t_{n}\rightarrow+\infty\}. \label{ForwardOmegaLlimitB}%
\end{equation}

\begin{definition}
\label{FAC}The process is forward asymptotically compact if for any sequence
$y_{n}\in U(t_{n},t_{0},B)$, where $t_{0}\in\mathbb{R}$ and $B$ is bounded,
there is a converging subsequence.
\end{definition}

\begin{lemma}
\label{ForwardOmega}If $U$ is forward asymptotically compact, then for any
bounded set $B$ and any $t_{0}\in\mathbb{R}$ the set $\omega_{f}(B,t_{0})$ is
non-empty, compact and forward attracts $B$, that is%
\begin{equation}
dist\left(  U(t,t_{0},B\right)  ,\omega_{f}(B,t_{0}))\rightarrow0\text{, as
}t\rightarrow+\infty\text{.} \label{ConvergOmega}%
\end{equation}
It is the minimal closed set satisfying (\ref{ConvergOmega}).
\end{lemma}

\begin{proof}
The fact that $\omega_{f}(B,t_{0})\not =\varnothing$ follows from Definition
\ref{FAC}.

Let $\{y_{m}\}\subset\omega_{f}(B,t_{0})$. Then there are $z_{m}\in
U(t_{m},t_{0},B),\ t_{m}\rightarrow+\infty$, such that $\rho(z_{m},y_{m}%
)\leq\frac{1}{m}$. Passing to a subsequence we have that $z_{m}\rightarrow
z_{0}\in\omega_{f}(B,t_{0})$. Hence,%
\[
\rho\left(  y_{m},z_{0}\right)  \leq\rho(y_{m},z_{m})+\rho(z_{m}%
,z_{0})\rightarrow0.
\]
Therefore, $\omega_{f}(B,t_{0})$ is a compact set.

Assume that (\ref{ConvergOmega}) is not true. Then there are $\varepsilon
>0,\ t_{n}\rightarrow+\infty$ and $y_{n}\in U(t_{n},t_{0},B)$ such that
\[
dist(y_{n},\omega_{f}(B,t_{0}))\geq\varepsilon\text{.}%
\]
However, there exists a subsequence $\{y_{n_{k}}\}$ such that $y_{n_{k}%
}\rightarrow y_{0}\in\omega_{f}(B,t_{0})$, which is a contradiction.

Finally, let $C$ be a closed set satisfying (\ref{ConvergOmega}). For
$y\in\omega_{f}(B,t_{0})$ let $y_{n}\in U(t_{n},t_{0},B)$ converges to $y$ as
$t_{n}\rightarrow+\infty$. Since $C$ forward attracts $B,$ for any
$\varepsilon>0$ there is $N_{1}(\varepsilon)$ such that%
\[
dist\left(  y_{n},C\right)  \leq\frac{\varepsilon}{2}\text{ if }n\geq N_{1}.
\]
Also, we choose $N_{2}\left(  \varepsilon\right)  $ such that $\rho\left(
y_{n},y\right)  \leq\frac{\varepsilon}{2}$ if $n\geq N_{2}$. Then%
\[
dist\left(  y,C\right)  \leq\rho\left(  y_{n},y\right)  +dist\left(
y,C\right)  \leq\varepsilon\text{ if }n\geq\max\{N_{1},N_{2}\}.
\]
Hence, $\omega_{f}(B,t_{0})\subset C.$
\end{proof}

\bigskip

For a family $\mathcal{K}=\{K(t)\}$ we define its $\omega$-limit set by%
\[
\omega\left(  \mathcal{K}\right)  =\{y:y=\lim_{n\rightarrow\infty}y_{n}\text{,
}y_{n}\in K(t_{n}),\ t_{n}\rightarrow+\infty\}.
\]
If the family $\mathcal{K}$ is invariant, then $K(t)=U(t,t_{0},K(t_{0}))$
implies that%
\[
\omega(\mathcal{K})=\omega_{f}(K(t_{0}),t_{0})\text{ for any }t_{0}%
\in\mathbb{R}\text{.}%
\]
Hence, if $U$ is forward asymptotically compact and $\mathcal{K}$ is bounded,
then Lemma \ref{ForwardOmega} implies that $\omega(\mathcal{K})$ is non-empty,
compact and
\begin{equation}
dist(K(t),\omega(\mathcal{K}))\rightarrow0\text{ as }t\rightarrow+\infty.
\label{ConvergOmega2}%
\end{equation}

We observe that%
\[
\omega\left(  \mathcal{K}\right)  =Limsup_{t\rightarrow+\infty}K(t)=\{y\in
X:\lim\inf_{t\rightarrow+\infty}dist(y,K(t))=0\}.
\]
We define also the set%
\begin{equation}
\omega_{0}\left(  \mathcal{K}\right)  =Liminf_{t\rightarrow+\infty}K(t)=\{y\in
X:\lim_{t\rightarrow+\infty}dist(y,K(t))=0\}. \label{DefOmega0}%
\end{equation}
These sets are closed. It is easy to see that $\omega_{0}(\mathcal{K}%
)\subset\omega(\mathcal{K})$. Also, if $\mathcal{K}$ is asymptotically
compact, then $\omega\left(  \mathcal{K}\right)  $ is non-empty while
$\omega_{0}\left(  \mathcal{K}\right)  $ can be empty in general.

\begin{lemma}
\label{ConvergOmega0}If for the family $\mathcal{K}$ the set $\omega
_{0}\left(  \mathcal{K}\right)  $ is non-empty and compact, then
$dist(\omega_{0}\left(  \mathcal{K}\right)  ,K(t))\rightarrow0$ as
$t\rightarrow+\infty. $
\end{lemma}

\begin{proof}
By contradiction assume that $dist(\omega_{0}\left(  \mathcal{K}\right)
,K(t))\not \rightarrow 0$ as $t\rightarrow+\infty.$ Then there are
$\varepsilon_{0}>0$, $t_{n}\rightarrow+\infty$ and $y_{n}\in\omega_{0}\left(
\mathcal{K}\right)  $ such that
\[
dist(y_{n},K(t_{n}))>\varepsilon_{0}\ \forall n.
\]
Since $\omega_{0}\left(  \mathcal{K}\right)  $ is compact, passing to a
subsequence we have that $y_{n}\rightarrow y_{0}\in\omega_{0}\left(
\mathcal{K}\right)  $. Then%
\[
dist(y_{n},K(t_{n}))\leq\rho(y_{n},y_{0})+dist(y_{0},K(t_{n}))\rightarrow
0\text{ as }n\rightarrow+\infty,
\]
which is a contradiction.
\end{proof}

\bigskip

We will give sufficient and necessary conditions for the existence of a
forward attractor in terms of the sets $\omega\left(  \mathcal{K}\right)  $
and $\omega_{0}\left(  \mathcal{K}\right)  .$

\begin{theorem}
\label{ForwardAttractorTh0}Let $U$ possess a pullback attractor $\mathcal{A}%
=\{\mathcal{A}(t)\}$. Assume that $U$ is forward asymptotically compact. If
for any bounded set $B$ and any $t_{0}\in\mathbb{R}$ we have
\begin{equation}
\omega_{f}(B,t_{0})\subset\omega_{0}(\mathcal{A}), \label{CondOmega0}%
\end{equation}
then $\mathcal{A}$ is a forward attractor and $\omega(\mathcal{A})=\omega
_{0}(\mathcal{A})$.

If $\mathcal{A}$ is a forward attractor, then%
\begin{equation}
\omega_{f}(B,t_{0})\subset\omega(\mathcal{A}),\text{ for any }t_{0}%
\in\mathbb{R}\text{ and }B\text{ bounded, } \label{CondOmega}%
\end{equation}

\end{theorem}

\begin{proof}
By Lemma \ref{ForwardOmega} and $U(t,t_{0},A(t_{0}))=A(t)$ we obtain that
$\omega(\mathcal{A})=\omega_{f}(A(t_{0}),t_{0})$ is non-empty and compact. If
(\ref{CondOmega0}) holds, then $\omega_{0}(\mathcal{A})$ is non-empty. Also,
$\omega_{0}(\mathcal{A})\subset\omega(\mathcal{A})$ implies that $\omega
_{0}(\mathcal{A})$ is compact. Therefore, by (\ref{CondOmega0}) and Lemmas
\ref{ForwardOmega} and \ref{ConvergOmega0} we deduce that
\begin{align*}
dist\left(  U(t,t_{0},B\right)  ,\mathcal{A}(t))  &  \leq dist(\left(
U(t,t_{0},B\right)  ,\omega_{0}(\mathcal{A}))+dist(\omega_{0}(\mathcal{A}%
),\mathcal{A}(t))\\
&  \leq dist(\left(  U(t,t_{0},B\right)  ,\omega_{f}(B,t_{0}))+dist(\omega
_{0}(\mathcal{A}),\mathcal{A}(t))\rightarrow0\text{, as }t\rightarrow+\infty.
\end{align*}
Hence, $\mathcal{A}$ is forward attracting. As $\omega_{0}(\mathcal{A}%
)\subset\omega(\mathcal{A})$, we need to establish the converse inclusion.
Using (\ref{CondOmega0}) we find
\[
\omega(\mathcal{A})=\omega_{f}(\mathcal{A}(t_{0}),t_{0})\subset\omega
_{0}(\mathcal{A}).
\]

Conversely, let $\mathcal{A}$ be a forward attractor. If (\ref{CondOmega})
were not true, there would exist $y_{0}\in\omega_{f}(B,t_{0})$ such that
$y_{0}\not \in \omega(\mathcal{A})$. Take a sequence $y_{n}\in U(t_{n}%
,t_{0},B)$, $t_{n}\rightarrow+\infty$, such that $y_{n}\rightarrow y_{0}$.
Since
\[
dist(y_{n},A(t_{n}))\rightarrow0,
\]
there are $z_{n}\in A(t_{n})$ satisfying that $\rho\left(  y_{n},z_{n}\right)
\rightarrow0$. This implies that $z_{n}\rightarrow y_{0}\in\omega
(\mathcal{A})$, which is a contradiction.
\end{proof}

\bigskip

We will show now with a simple example that condition (\ref{CondOmega0}) is
not necessary for the existence of the forward attractor. We consider the
equation%
\[
\frac{dy}{dt}=y+\sin(t),
\]
whose solution is
\[
y(t)=e^{t_{0}-t}y_{0}+\int_{t_{0}}^{t}e^{s-t}\sin(s)ds.
\]
The pullback attractor is given by%
\[
A(t)=\{a(t)\},
\]%
\[
a(t)=\int_{-\infty}^{t}e^{s-t}\sin(s)ds=\frac{1}{2}\left(  \sin(t\right)
-\cos(t)).
\]
It is easy to check that it is also a forward attractor and that%
\begin{align*}
\omega_{0}(\mathcal{A})  &  =\varnothing,\\
\omega(\mathcal{A})  &  =[-b,b],\ b=\frac{1}{2}\left(  \sin(\frac{3\pi}%
{4})-\cos(\frac{3\pi}{4})\right)  =\frac{\sqrt{2}}{2}.
\end{align*}
Thus, (\ref{CondOmega0}) is not satisfied.

\begin{lemma}
\label{EquivConditions}Let $U$ possess a pullback attractor $\mathcal{A}%
=\{\mathcal{A}(t)\}$. Assume that $U$ is forward asymptotically compact. Then
condition (\ref{CondOmega0}) is satisfied if and only if (\ref{CondOmega}) and
the condition%
\begin{equation}
dist_{H}(\mathcal{A}(t),\omega(\mathcal{A}))\rightarrow0,\text{ as
}t\rightarrow+\infty, \label{CondOmega2}%
\end{equation}
hold true.
\end{lemma}

\begin{proof}
Let (\ref{CondOmega0}) be satisfied. Then $\omega_{0}(\mathcal{A}%
)\subset\omega(\mathcal{A})$ implies (\ref{CondOmega}). The convergence
$dist(\mathcal{A}(t),\omega(\mathcal{A}))\rightarrow0,$ as $t\rightarrow
+\infty,$ is a consequence of (\ref{ConvergOmega2}). From Theorem
\ref{ForwardAttractorTh0} we have that $\omega(\mathcal{A})\subset\omega
_{0}(\mathcal{A})$, and, therefore, Lemma \ref{ConvergOmega0} yields%
\begin{equation}
dist(\omega(\mathcal{A}),\mathcal{A}(t))\rightarrow0\text{ as }t\rightarrow
+\infty. \label{CondOmega3}%
\end{equation}

Conversely, let (\ref{CondOmega}) and (\ref{CondOmega2}) be satisfied. From
(\ref{CondOmega2}) and the definition of $\omega_{0}(\mathcal{A})$ we obtain
that $\omega(\mathcal{A})\subset\omega_{0}(\mathcal{A})$, Hence,
$\omega(\mathcal{A})=\omega_{0}(\mathcal{A})$ and (\ref{CondOmega}) give
(\ref{CondOmega0}).
\end{proof}

\begin{theorem}
\label{ForwardAttractorTh}Let $U$ possess a pullback attractor $\mathcal{A}%
=\{\mathcal{A}(t)\}$. Assume that $U$ is forward asymptotically compact and
that (\ref{CondOmega}), (\ref{CondOmega2}) are satisfied. Then $\mathcal{A}$
is a forward attractor.
\end{theorem}

\begin{proof}
It follows from Lemma \ref{EquivConditions} and Theorem
\ref{ForwardAttractorTh0}.
\end{proof}

\begin{remark}
The conditions in Theorems \ref{ForwardAttractorTh0}, \ref{ForwardAttractorTh}
are strong in the sense that they imply the equality $\omega_{0}%
(\mathcal{A})=\omega(\mathcal{A})$. However, we have seen in the previous
example that a forward attractor can exist when the strict inclusion
$\omega_{0}(\mathcal{A})\subset\omega(\mathcal{A})$ holds.
\end{remark}

Following \cite{CuiFigLangaNasc} we define the set%
\begin{equation}
\mathcal{A}_{\min}=\overline{\underset{B\text{ bounded}}{\cup}\underset{t_{0}%
\in\mathbb{R}}{\cup}\omega_{f}(B,t_{0})}. \label{Amin}%
\end{equation}
Lemma \ref{ForwardOmega} implies that $\mathcal{A}_{\min}$ is the minimal
closed set forward attracting any bounded set $B$. The condition%
\begin{equation}
dist\left(  \mathcal{A}_{\min},\mathcal{A}(t)\right)  \rightarrow0\text{, as
}t\rightarrow+\infty, \label{CondAmin}%
\end{equation}
was used in \cite{CuiFigLangaNasc} to establish that the family $\mathcal{A}$
is a forward attractor.

\begin{lemma}
\label{EquivSuffCondition}Let $U$ possess a pullback attractor $\mathcal{A}%
=\{\mathcal{A}(t)\}$. Assume that $U$ is forward asymptotically compact. Then
condition (\ref{CondOmega0}) is satisfied if and only (\ref{CondAmin}) holds true.
\end{lemma}

\begin{proof}
Let (\ref{CondOmega0}) be satisfied. By Theorem \ref{ForwardAttractorTh0} we
deduce that $\omega_{0}(\mathcal{A})$ is non-empty and compact. Then
$\mathcal{A}_{\min}\subset\omega_{0}(\mathcal{A})$ implies by Lemma
\ref{ConvergOmega0} that%
\begin{align*}
dist\left(  \mathcal{A}_{\min},\mathcal{A}(t)\right)   &  \leq dist\left(
\mathcal{A}_{\min},\omega_{0}(\mathcal{A})\right)  +dist\left(  \omega
_{0}(\mathcal{A}),\mathcal{A}(t)\right) \\
&  =dist\left(  \omega_{0}(\mathcal{A}),\mathcal{A}(t)\right)  \rightarrow
0\text{, as }t\rightarrow+\infty.
\end{align*}
Conversely, let (\ref{CondAmin}) hold true. Then any $y\in\omega_{f}(B,t_{0})$
satisfies%
\[
dist\left(  y,\mathcal{A}(t)\right)  \rightarrow0\text{, as }t\rightarrow
+\infty.
\]
By the definition (\ref{DefOmega0}) we derive that $y\in\omega_{0}%
(\mathcal{A})$. Therefore, (\ref{CondOmega0}) holds true.
\end{proof}

\bigskip

We observe also that since $\omega\left(  \mathcal{A}\right)  =\omega
_{f}\left(  \mathcal{A}\left(  t_{0}\right)  ,t_{0}\right)  \subset
\mathcal{A}_{\min}$, the necessary condition (\ref{CondOmega}) can be written
in the following equivalent way:%
\[
\mathcal{A}_{\min}=\omega\left(  \mathcal{A}\right)  .
\]
This is the necessary condition given in \cite[Theorem 6.4]{CuiFigLangaNasc}.

\bigskip

Let us consider now the specific case of asymptotically autonomous systems.

We consider a continuous semigroup $S:\mathbb{R}^{+}\times X\rightarrow X$,
which satisfies:

\begin{enumerate}
\item $S(0,x)=x$ for all $x\in X,;$

\item $S(t+s,x)=S(t,S(s,x))$ for all $x\in X$, $0\leq s\leq t;$

\item $\left(  t,x\right)  \mapsto S(t,x)$ is a continuous map.
\end{enumerate}

A compact set $\mathcal{A}_{\infty}$ is said to be a global attractor for $S $
if it is invariant, that is, $S(t,\mathcal{A}_{\infty})=\mathcal{A}_{\infty}$,
for any $t\geq0$, and attracts any bounded set $B$, which means that%
\[
dist(S(t,B),\mathcal{A}_{\infty})\rightarrow0\text{ as }t\rightarrow+\infty.
\]
A complete trajectory of $S$ is a function $\varphi:\mathbb{R}\rightarrow X$
such that $\varphi\left(  t\right)  =S(t-s,\varphi(s))$ for all $t\geq s\geq
0$. A complete trajectory $\varphi$ is bounded if $\cup_{t\in\mathbb{R}%
}\varphi\left(  t\right)  $ is a bounded set. The global attractor is
characterized by the union of all bounded complete trajectories \cite[p.10]%
{LadBook}, that is,%
\begin{equation}
\mathcal{A}_{\infty}=\{\varphi\left(  0\right)  :\varphi\text{ is a bounded
complete trajectory}\}. \label{CharacAttr}%
\end{equation}

We will say that the process is asymptotically autonomous if there exists a
continuous semigroup $S:\mathbb{R}^{+}\times X\rightarrow X$ such that for any
sequences $\tau_{n}\rightarrow+\infty,\ x_{n}\rightarrow x_{0}$, as
$n\rightarrow+\infty,$ we have that%
\begin{equation}
U(t+\tau_{n},\tau_{n},x_{n})\rightarrow S(t,x_{0})\text{, as }n\rightarrow
+\infty\text{, for all }t\geq0. \label{AA}%
\end{equation}

\begin{lemma}
\label{ConvergAuton}Let the process $U$ have the pullback attractor
$\mathcal{A}=\{\mathcal{A}(t)\}$ and let the semigroup $S:\mathbb{R}^{+}\times
X\rightarrow X$ possess the global attractor $\mathcal{A}_{\infty}$. Let $U$
be asymptotically autonomous and forward asymptotically compact. Then%
\[
dist(\mathcal{A}(t),\mathcal{A}_{\infty})\rightarrow0\ \text{as }%
t\rightarrow+\infty.
\]

\end{lemma}

\begin{proof}
Assume that it is not true. Then there would exist $\varepsilon_{0}>0$,
$t_{n}\nearrow+\infty$ and $y_{n}\in\mathcal{A}(t_{n})$ such that%
\[
dist(y_{n},\mathcal{A}_{\infty})\geq\varepsilon_{0}\text{ for all }n.
\]
By (\ref{ConvergOmega2}) we know that $dist(\mathcal{A}(t),\omega
(\mathcal{A}))\rightarrow0$. Therefore, there is $t_{0}$ such that
$C_{0}=\overline{\cup_{t\geq t_{0}}\mathcal{A}(t)}$ is bounded. The attraction
property of $\mathcal{A}_{\infty}$ implies the existence of $n_{0}\left(
\varepsilon_{0}\right)  $ such that%
\[
dist(S(t_{n},C_{0}),\mathcal{A}_{\infty})\leq\frac{\varepsilon_{0}}{3}\text{
if }n\geq n_{0}.
\]
Since $U(t_{n},t_{n}-t_{n_{0}},\mathcal{A}(t_{n}-t_{n_{0}}))=\mathcal{A}%
(t_{n})$ for $n\geq n_{0}$, there are $z_{n}\in\mathcal{A}(t_{n}-t_{n_{0}})$
such that $y_{n}=U(t_{n},t_{n}-t_{n_{0}},z_{n})$ and, up to a subsequence,
$z_{n}\rightarrow z_{0}\in C_{0}$. Let $\tau_{n}=t_{n}-t_{n_{0}}$ By
(\ref{AA}) there is $n_{1}\left(  \varepsilon_{0}\right)  \geq n_{0}$ such
that%
\[
\rho(y_{n},S(t_{n_{0}},z_{0}))=\rho(U(\tau_{n}+t_{n_{0}},\tau_{n}%
,z_{n}),S(t_{n_{0}},z_{0}))\leq\frac{\varepsilon_{0}}{3}\text{ for }n\geq
n_{1}.
\]
Thus,%
\[
dist(y_{n},\mathcal{A}_{\infty})\leq dist(y_{n},S(t_{n_{0}},z_{0}%
))+\rho(S(t_{n_{0}},z_{0}),\mathcal{A}_{\infty})\leq\frac{2\varepsilon_{0}}%
{3},\text{ for }n\geq n_{1},
\]
which is a contradiction.
\end{proof}

\begin{remark}
This lemma is a version of Theorem 4.1 in \cite{KloedenLorenz} with slightly
different conditions.
\end{remark}

\begin{theorem}
\label{ForwardAttractorTh2}Let the process $U$ have the pullback attractor
$\mathcal{A}=\{\mathcal{A}(t)\}$ and let the semigroup $S:\mathbb{R}^{+}\times
X\rightarrow X$ possess the global attractor $\mathcal{A}_{\infty}$. Let $U$
be asymptotically autonomous and forward asymptotically compact. Assume that%
\begin{equation}
dist_{H}(\mathcal{A}(t),\mathcal{A}_{\infty})\rightarrow0\text{ as
}t\rightarrow+\infty. \label{CondAInf}%
\end{equation}
Then $\mathcal{A}$ is a forward attractor.
\end{theorem}

\begin{proof}
We will check that conditions (\ref{CondOmega})-(\ref{CondOmega2}) hold, and
then the result follows from Theorem \ref{ForwardAttractorTh}.

In view of (\ref{ConvergOmega2}), for assumption (\ref{CondOmega2}) we just
need to prove (\ref{CondOmega3}). We observe that $\omega(\mathcal{A}%
)\subset\mathcal{A}_{\infty}$. Indeed, if $y\in\omega(\mathcal{A})$, there are
$y_{n}\in\mathcal{A}(t_{n})$, $t_{n}\rightarrow+\infty$, such that
$y_{n}\rightarrow y$. From (\ref{CondAInf}) we have $dist(y_{n},\mathcal{A}%
_{\infty})\rightarrow0$ and, therefore, the compactness of $\mathcal{A}%
_{\infty}$ implies that $y\in\mathcal{A}_{\infty}$. Then, by (\ref{CondAInf})
we infer%
\[
dist(\omega(\mathcal{A}),\mathcal{A}(t))\leq dist(\mathcal{A}_{\infty
},\mathcal{A}(t))\rightarrow0.
\]

Let us prove (\ref{CondOmega}). As conditions (\ref{CondOmega2}) and
(\ref{CondAInf}) are satisfied, it is easy to see that $\omega(\mathcal{A}%
)=\mathcal{A}_{\infty}$. If $y_{0}\in\omega_{f}(B,t_{0})$, there are
$y_{n}=u_{n}(t_{n})=U(t_{n},t_{0},x_{n})$, with $t_{n}\rightarrow
+\infty,\ x_{n}\in B$, such that $y_{n}\rightarrow y_{0}$. Let $v_{n}%
(t)=u_{n}(t+t_{n})$. Then $v_{n}(0)\rightarrow y_{0}$ and
\[
v_{n}(t)=U(t_{n}+t,t_{0},x_{n})=U(t_{n}+t,t_{n},U(t_{n},t_{0},x_{n}%
))=U(t_{n}+t,t_{n},y_{n})\rightarrow v_{0}(t)\text{ for any }t\geq0,
\]
where $v_{0}(t)=S(t,y_{0})\in\omega_{f}(B,t_{0}).$ Put $\psi_{0}(t)=v_{0}(t)$
for $t\geq0$. Further, let $v_{n}^{-1}(t)=u_{n}(t+t_{n}-1)=U(t_{n}%
+t-1,t_{0},x_{n})$. Then, since $U$ is forward asymptotically compact, up to a
subsequence $v_{n}^{-1}(0)\rightarrow y^{-1}\in\omega_{f}(B,t_{0})$ and%
\[
v_{n}^{-1}(t)=U(t_{n}+t-1,t_{0},x_{n})=U(t_{n}+t-1,t_{n}-1,U(t_{n}%
-1,t_{0},x_{n}))=U(t_{n}+t-1,t_{n}-1,v_{n}^{-1}(0))\rightarrow v_{-1}%
(t)\text{,}%
\]
for any $t\geq0,$ where $v_{-1}(t)=S(t,y_{-1})\in\omega_{f}(B,t_{0}).$ We put
$\psi_{-1}(t)=v_{-1}(t+1)$. Then $\psi_{-1}(t)=\psi_{0}(t),$ for any $t\geq0,$
and
\[
S(t-s,\psi_{-1}(s))=S(t-s,v_{-1}(s+1))
\]%
\[
=S(t-s,S(s+1,y_{-1}))=S(t+1,y_{-1})=v_{-1}(t+1)=\psi_{-1}(t),\text{ for
all\ }t\geq s\geq-1.
\]
Proceeding in this way for $k=-2,-3,...$ we obtain a sequence of functions
$\psi_{-k}:[-k,\infty)\rightarrow X$, $k=0,1,2,...,$ such that $\psi
_{-k}(t)\in\omega_{f}(B,t_{0})$, for $t\geq-k$, $\psi_{-k}(t)=\psi_{-k+1}(t),$
for any $t\geq-k+1,$ and
\[
S(t-s,\psi_{-k}(s))=S(t-s,v_{-k}(s+k))
\]%
\[
=S(t-s,S(s+k,y_{-k}))=S(t+k,y_{-k})=v_{-k}(t+k)=\psi_{-k}(t),\text{ for
all\ }t\geq s\geq-k.
\]
We define $\psi:\left(  -\infty,+\infty\right)  \rightarrow X$ by $\psi\left(
t\right)  =\psi_{-k}(t)$, for $t\geq-k$, for any $k=0,1,2,...$ Then $\psi$ is
a bounded complete trajectory of the semigroup $S$ and $\psi\left(  0\right)
=y_{0}$. The characterization (\ref{CharacAttr}) of the global attractor
implies then that $y_{0}\in\mathcal{A}_{\infty}=\omega(\mathcal{A})$. Thus,
$\omega_{f}(B,t_{0})\subset\omega(\mathcal{A})$.
\end{proof}

\bigskip

\subsection{Multivalued processes}

Let us denote by $P(X)$ ($C(X)$) the set of all non-empty (non-empty closed)
subsets of $X$. We consider a multivalued process $U:\mathbb{R}_{\geq}%
^{2}\times X\rightarrow P(X)$, which satisfies:

\begin{enumerate}
\item $U(s,s,x)=x$ for all $x\in X,\ s\in\mathbb{R};$

\item $U(t,s,x)\subset U(t,\tau,U(\tau,s,x))$ for all $x\in X,\ t,s,\tau
\in\mathbb{R}$, $s\leq\tau\leq t.$
\end{enumerate}

If, in addition, $U(t,s,x)=U(t,\tau,U(\tau,s,x))$ for all $x\in X,\ t,s,\tau
\in\mathbb{R}$, $s\leq\tau\leq t$, then the process is said to be strict.

A family $\{K(t)\}$ is negatively (positively) invariant if $K(t)\subset
U(t,s,K(s))$ ($K(t)\supset U(t,s,K(s))$) for any $\left(  t,s\right)
\in\mathbb{R}_{\geq}^{2}$. It is invariant if $K(t)=U(t,s,K(s))$ for any
$\left(  t,s\right)  \in\mathbb{R}_{\geq}^{2}$.

\begin{definition}
A family of non-empty sets $\mathcal{A}=\{\mathcal{A}(t)\}$ is called pullback
attracting if it attracts any bounded set $B\subset X$ in the pullback sense,
that is,%
\[
dist\left(  U(t,t_{0},B\right)  ,\mathcal{A}(t))\rightarrow0\text{, as }%
t_{0}\rightarrow-\infty\text{, for any }t\text{ fixed.}%
\]
It is called forward attracting if it attracts any bounded set $B\subset X$ in
the forward sense, that is,%
\[
dist\left(  U(t,t_{0},B\right)  ,\mathcal{A}(t))\rightarrow0\text{, as
}t\rightarrow+\infty\text{, for any }t_{0}\text{ fixed.}%
\]

\end{definition}

\begin{definition}
A family of non-empty compact negatively invariant sets $\mathcal{A}%
=\{\mathcal{A}(t)\}$ is called a pullback attractor if it is the minimal
pullback attracting family of closed sets.
\end{definition}

\begin{definition}
A family of non-empty compact negatively invariant sets $\mathcal{A}%
=\{\mathcal{A}(t)\}$ is called a forward attractor if it is forward attracting.
\end{definition}

As in the single-valued situation, although several forward attractors could
exist, they are asymptotically equivalent.

\begin{lemma}
If $\mathcal{A}_{1},\mathcal{A}_{2}$ are two forward attractors, then
\[
dist_{H}(\mathcal{A}_{1}\left(  t\right)  ,\mathcal{A}_{2}\left(  t\right)
)\rightarrow0\text{ as }t\rightarrow+\infty.
\]

\end{lemma}

\begin{proof}
Take $\tau\in\mathbb{R}$. Since $\mathcal{A}_{1}\left(  t\right)  \subset
U(t,\tau,\mathcal{A}_{1}\left(  \tau\right)  ),$ we have that
\[
dist\left(  \mathcal{A}_{1}\left(  t\right)  ,\mathcal{A}_{2}\left(  t\right)
\right)  \leq dist\left(  U(t,\tau,\mathcal{A}_{1}\left(  \tau\right)
),\mathcal{A}_{2}\left(  \tau\right)  \right)  \rightarrow0\text{, as
}t\rightarrow+\infty.
\]
Likewise we obtain that $dist\left(  \mathcal{A}_{2}\left(  t\right)
,\mathcal{A}_{1}\left(  t\right)  \right)  \rightarrow0$ as $t\rightarrow
+\infty.$
\end{proof}

\bigskip

General results for existence of pullback attractors for multivalued processes
can be found in \cite{CLMV03}, \cite{CotiZelatiKalita}, \cite{KapKasVal2011},
while their characterization was studied in \cite{CLV16}. If the process $U$
is strict and the pullback attractor $\mathcal{A}$ is backward bounded, that
is, $\cup_{t\leq\tau}\mathcal{A}(t)$ is bounded for some $\tau\in\mathbb{R}$,
then it is known \cite[Lemma 2.5]{CLV16} that it is invariant.

The concept of forward $\omega$-limit set and Definition \ref{FAC} are the
same as in the single-valued case.

\begin{lemma}
\label{ForwardOmegaM}If $U$ is forward asymptotically compact, then for any
bounded set $B$ and any $t_{0}\in\mathbb{R}$ the set $\omega_{f}(B,t_{0})$ is
non-empty, compact and forward attracts $B$, that is,%
\begin{equation}
dist\left(  U(t,t_{0},B\right)  ,\omega_{f}(B,t_{0}))\rightarrow0\text{, as
}t\rightarrow+\infty\text{.} \label{ConvergOmegaM}%
\end{equation}
It is the minimal closed set satisfying (\ref{ConvergOmegaM}).
\end{lemma}

\begin{proof}
It is the same as in Lemma \ref{ForwardOmega}.
\end{proof}

\bigskip

\begin{lemma}
\label{ForwardA}If the family $\mathcal{K}=\{K(t)\}$ is forward asymptotically
compact, then the set $\omega(\mathcal{K})$ is non-empty, compact and forward
attracts $\mathcal{K}$, that is%
\[
dist(K(t),\omega(\mathcal{K}))\rightarrow0\text{, as }t\rightarrow
+\infty\text{.}%
\]

\end{lemma}

\begin{proof}
The proof is quite similar to the one in Lemma \ref{ForwardOmega}.
\end{proof}

\bigskip

As in the single-valued case, we will give sufficient and necessary conditions
for the existence of a forward attractor in terms of the sets $\omega\left(
\mathcal{A}\right)  $ and $\omega_{0}\left(  \mathcal{A}\right)  .$

\begin{theorem}
\label{ForwardAttractorTh0M}Let $U$ possess a pullback attractor
$\mathcal{A}=\{\mathcal{A}(t)\}$. Assume that $U$ is forward asymptotically
compact. If for any bounded set $B$ and any $t_{0}\in\mathbb{R}$ we have
\begin{equation}
\omega_{f}(B,t_{0})\subset\omega_{0}(\mathcal{A}), \label{CondOmega0M}%
\end{equation}
then $\mathcal{A}$ is a forward attractor and $\omega(\mathcal{A})=\omega
_{0}(\mathcal{A})$.

If $\mathcal{A}$ is a forward attractor, then%
\begin{equation}
\omega_{f}(B,t_{0})\subset\omega(\mathcal{A}),\text{ for any }t_{0}%
\in\mathbb{R}\text{ and }B\text{ bounded, } \label{CondOmegaM}%
\end{equation}

\end{theorem}

\begin{proof}
The inclusion $\mathcal{A}(t)\subset U(t,t_{0},\mathcal{A}(t_{0}))$ implies
that $\mathcal{A}$ is asymptotically compact. By Lemma \ref{ForwardA} we
obtain that $\omega(\mathcal{A})$ is non-empty and compact and condition
(\ref{CondOmega0M}) implies that $\omega_{0}(\mathcal{A})$ is non-empty. By
$\omega_{0}(\mathcal{A})\subset\omega(\mathcal{A})$ we find that $\omega
_{0}(\mathcal{A})$ is compact. Therefore, by (\ref{CondOmega0M}) and Lemmas
\ref{ForwardOmegaM} and \ref{ConvergOmega0} we deduce that
\begin{align*}
dist\left(  U(t,t_{0},B\right)  ,\mathcal{A}(t))  &  \leq dist(\left(
U(t,t_{0},B\right)  ,\omega_{0}(\mathcal{A}))+dist(\omega_{0}(\mathcal{A}%
),\mathcal{A}(t))\\
&  \leq dist(\left(  U(t,t_{0},B\right)  ,\omega_{f}(B,t_{0}))+dist(\omega
_{0}(\mathcal{A}),\mathcal{A}(t))\rightarrow0\text{, as }t\rightarrow+\infty.
\end{align*}

As $\omega_{0}(\mathcal{A})\subset\omega(\mathcal{A})$, we only need to prove
that $\omega(\mathcal{A})\subset\omega_{0}(\mathcal{A})$. Since $\mathcal{A}%
(t)\subset U(t,t_{0},\mathcal{A}(t_{0}))$, by (\ref{CondOmega0M}) we obtain
\[
\omega(\mathcal{A})\subset\omega_{f}(\mathcal{A}(t_{0}),t_{0})\subset
\omega_{0}(\mathcal{A}).
\]

Let now $\mathcal{A}$ be a forward attractor. Let there be $y_{0}\in\omega
_{f}(B,t_{0})$ such that $y_{0}\not \in \omega(\mathcal{A})$. Take a sequence
$y_{n}\in U(t_{n},t_{0},B)$, $t_{n}\rightarrow+\infty$, such that
$y_{n}\rightarrow y_{0}$. Since
\[
dist(y_{n},\mathcal{A}(t_{n}))\rightarrow0,
\]
there are $z_{n}\in\mathcal{A}(t_{n})$ satisfying that $\rho\left(
y_{n},z_{n}\right)  \rightarrow0$. Hence, $z_{n}\rightarrow y_{0}\in
\omega(\mathcal{A})$, which is a contradiction.
\end{proof}

\begin{lemma}
\label{EquivConditionsM}Let $U$ possess a pullback attractor $\mathcal{A}%
=\{\mathcal{A}(t)\}$. Assume that $U$ is forward asymptotically compact. Then
condition (\ref{CondOmega0M}) is satisfied if and only if (\ref{CondOmegaM})
and
\begin{equation}
dist_{H}(\mathcal{A}(t),\omega(\mathcal{A}))\rightarrow0,\text{ as
}t\rightarrow+\infty, \label{CondOmega2M}%
\end{equation}
hold true.
\end{lemma}

\begin{proof}
If (\ref{CondOmega0M}) is satisfied, then $\omega_{0}(\mathcal{A}%
)\subset\omega(\mathcal{A})$ gives (\ref{CondOmegaM}). The inclusion
$\mathcal{A}(t)\subset U(t,t_{0},\mathcal{A}(t_{0}))$ implies that
$\mathcal{A}$ is asymptotically compact and, therefore, Lemma \ref{ForwardA}
implies that $dist(\mathcal{A}(t),\omega(\mathcal{A}))\rightarrow0,$ as
$t\rightarrow+\infty$. From Theorem \ref{ForwardAttractorTh0M} we have that
$\omega(\mathcal{A})\subset\omega_{0}(\mathcal{A})$, and, therefore, from
Lemma \ref{ConvergOmega0} we have (\ref{CondOmega2M}).

Conversely, let (\ref{CondOmegaM}) and (\ref{CondOmega2M}) hold. From
(\ref{CondOmega2M}) and the definition of $\omega_{0}(\mathcal{A})$ we deduce
that $\omega(\mathcal{A})\subset\omega_{0}(\mathcal{A})$, Hence,
$\omega(\mathcal{A})=\omega_{0}(\mathcal{A})$ and in such a case
(\ref{CondOmegaM}) and (\ref{CondOmega0M}) are equivalent.
\end{proof}

\begin{theorem}
\label{ForwardAttractorThM}Let $U$ possess a pullback attractor $\mathcal{A}%
=\{\mathcal{A}(t)\}$. Assume that $U$ is forward asymptotically compact and
that (\ref{CondOmegaM}), (\ref{CondOmega2M}) are satisfied. Then $\mathcal{A}$
is a forward attractor.
\end{theorem}

\begin{proof}
It follows from Lemma \ref{EquivConditionsM} and Theorem
\ref{ForwardAttractorTh0M}.
\end{proof}

\bigskip

As in the single-valued situation we define the set (\ref{Amin}), which is by
Lemma \ref{ForwardOmegaM} the minimal closed set forward attracting any
bounded set $B$. We consider then the condition%
\begin{equation}
dist\left(  \mathcal{A}_{\min},\mathcal{A}(t)\right)  \rightarrow0\text{, as
}t\rightarrow+\infty. \label{CondAminM}%
\end{equation}

\begin{lemma}
Let $U$ possess a pullback attractor $\mathcal{A}=\{\mathcal{A}(t)\}$. Assume
that $U$ is forward asymptotically compact. Then condition (\ref{CondOmega0M})
is satisfied if and only (\ref{CondAminM}) holds true.
\end{lemma}

\begin{proof}
It repeats the same steps of the proof of Lemma \ref{EquivSuffCondition} but
using Theorem \ref{ForwardAttractorTh0M}.
\end{proof}

\bigskip

In order to obtain the equivalent result from Theorem
\ref{ForwardAttractorTh2}, as we need to work with the structure of the
attractor, we will define the process and the semigroup in terms of trajectories.

Let $W_{\tau}=C([\tau,\infty),X)$ and consider a family of functions (called
trajectories) $\mathcal{R}=\{\mathcal{R}_{\tau}\}_{\tau\in\mathbb{R}}$, where
$\mathcal{R}_{\tau}\subset W_{\tau}$. Let us define the following axioms:

\begin{itemize}
\item[$\left(  K1\right)  $] (Existence property) For any $\tau\in
\mathbb{R},\ x\in X$ there exists $\varphi\in\mathcal{R}_{\tau}$ satisfying
$\varphi\left(  \tau\right)  =x.$

\item[$\left(  K2\right)  $] (Translation property) For any $s>0,\ \varphi
\in\mathcal{R}_{\tau}$, the function $\varphi_{s}=\varphi\mid_{\lbrack
\tau+s,\infty)}$ belongs to $\mathcal{R}_{\tau+s}.$

\item[$\left(  K3\right)  $] (Concatenation property) If $\varphi_{1}%
\in\mathcal{R}_{\tau},\varphi_{2}\in\mathcal{R}_{s}$, where $s>\tau$, are such
that $\varphi_{1}\left(  s\right)  =\varphi_{2}(s)$, then%
\[
\varphi\left(  t\right)  =\left\{
\begin{array}
[c]{c}%
\varphi_{1}(t)\text{ if }t\in\lbrack\tau,s],\\
\varphi_{2}(t)\text{ if }t\geq s,
\end{array}
\right.
\]
belongs to $\mathcal{R}_{\tau}$.

\item[$\left(  K4\right)  $] (Continuity property) If $\{\varphi_{n}%
\}\subset\mathcal{R}_{\tau}$ and $\varphi_{n}(\tau)\rightarrow\varphi_{0}$,
then there is a subsequence $\{\varphi_{n_{k}}\}$ and $\varphi\in
\mathcal{R}_{\tau}$ such that $\varphi\left(  \tau\right)  =\varphi_{0}$ and%
\[
\varphi_{n_{k}}\left(  t\right)  \rightarrow\varphi\left(  t\right)  \text{
for any }t\geq\tau.
\]

\end{itemize}

Properties (K1)-(K2) imply that the multivalued map $U:\mathbb{R}_{\geq}%
^{2}\times X\rightarrow P(X)$ defined by%
\[
U(t,t_{0},x)=\{y:y=\varphi\left(  t\right)  \text{ for some }\varphi
\in\mathcal{R}_{\tau}\text{ such that }\varphi\left(  \tau\right)  =x\}
\]
is a multivalued process. If $\left(  K1\right)  -\left(  K3\right)  $ hold,
then $U$ is a strict multivalued process.

A complete trajectory of $\mathcal{R}$ is a function $\xi:\mathbb{R}%
\rightarrow X$ such that $\xi\mid_{\lbrack s,\infty)}\in\mathcal{R}_{\tau}$
for all $s\in\mathbb{R}.$ A complete trajectory $\xi$ is bounded if
$\cup_{t\in\mathbb{R}}\xi\left(  t\right)  $ is a bounded set. If $\left(
H1\right)  -\left(  H2\right)  $ are satisfied, either $\left(  H3\right)  $
or $\left(  H4\right)  $ holds, and the pullback attractor $\mathcal{A}$ is
globally bounded, that is, $\cup_{t\in\mathbb{R}}\mathcal{A}(t)$ is bounded,
then $\mathcal{A}$ can be characterized as the union of of bounded complete
trajectories \cite[Corollaries 2.10 and 2.12 ]{CLV16}, that is,%
\[
\mathcal{A}(t)=\{\xi\left(  t\right)  :\xi\text{ is a bounded complete
trajectory}\}.
\]

In the same way, for the autonomous case let us consider a family of functions
$\mathcal{R}_{0}\subset W_{0}$ and the axioms:

\begin{itemize}
\item[$\left(  H1\right)  $] (Existence property) For any $\ x\in X$ there
exists $\varphi\in\mathcal{R}_{0}$ satisfying $\varphi\left(  0\right)  =x.$

\item[$\left(  H2\right)  $] (Translation property) For any $s>0,\ \varphi
\in\mathcal{R}_{0}$, the function $\varphi_{s}\left(
\text{\textperiodcentered}\right)  =\varphi($\textperiodcentered$+s)$ belongs
to $\mathcal{R}_{0}.$

\item[$\left(  H3\right)  $] (Concatenation property) If $\varphi_{1}%
\in\mathcal{R}_{0},\varphi_{2}\in\mathcal{R}_{0}$, where $s>0$, are such that
$\varphi_{1}\left(  s\right)  =\varphi_{2}(0)$, then%
\[
\varphi\left(  t\right)  =\left\{
\begin{array}
[c]{c}%
\varphi_{1}(t)\text{ if }t\in\lbrack0,s],\\
\varphi_{2}(t-s)\text{ if }t\geq s,
\end{array}
\right.
\]
belongs to $\mathcal{R}_{0}$.

\item[$\left(  H4\right)  $] (Continuity property) If $\{\varphi_{n}%
\}\subset\mathcal{R}_{0}$ and $\varphi_{n}(0)\rightarrow\varphi_{0}$, then
there is a subsequence $\{\varphi_{n_{k}}\}$ and $\varphi\in\mathcal{R}_{0}$
such that $\varphi\left(  0\right)  =\varphi_{0}$ and%
\[
\varphi_{n_{k}}\left(  t\right)  \rightarrow\varphi\left(  t\right)  \text{
for any }t\geq0.
\]

\end{itemize}

We recall that $G:\mathbb{R}^{+}\times X\rightarrow P(X)$ is a multivalued
semiflow if:

\begin{itemize}
\item $G(0,$\textperiodcentered$)$ is the identity map;

\item $G(t+s,x)\subset G(t,G(s,x))$ for all $x\in X,\ 0\leq s\leq t.$
\end{itemize}

If, moreover, $G(t+s,x)=G(t,G(s,x))$, then $G$ is a strict semiflow. Axioms
$\left(  H1\right)  -\left(  H2\right)  $ imply that the multivalued map
$G:\mathbb{R}^{+}\times X\rightarrow P(X)$ defined by%
\[
G(t,x)=\{y:y=\varphi\left(  t\right)  \text{ for some }\varphi\in
\mathcal{R}_{0}\text{ such that }\varphi\left(  0\right)  =x\}
\]
is a multivalued semiflow. If $\left(  H1\right)  -\left(  H3\right)  $ hold,
then $G$ is a strict multivalued semiflow.

A compact set $\mathcal{A}_{\infty}$ is said to be a global attractor for $G$
if it is negatively invariant, that is, $\mathcal{A}_{\infty}\subset
G(t,\mathcal{A}_{\infty})$, for all $t\geq0$, and attracts any bounded set
$B$, that is,%
\[
dist\left(  G(t,B),\mathcal{A}_{\infty}\right)  \rightarrow0\text{ as
}t\rightarrow+\infty.
\]
If $G$ is a strict semiflow, then the global attractor $\mathcal{A}_{\infty}$
is invariant, i.e., $\mathcal{A}_{\infty}=G(t,\mathcal{A}_{\infty})$, for all
$t\geq0$ \cite[Remark 8]{MV98}.

A function $\phi:\mathbb{R}\rightarrow X$ is a complete trajectory of
$\mathcal{R}_{0}$ if $\phi\left(  \text{\textperiodcentered}+s\right)
\in\mathcal{R}_{0}$ for any $s\in\mathbb{R}.$ It is said to be bounded if
$\cup_{t\in\mathbb{R}}\phi\left(  t\right)  $ is a bounded set. If $\left(
H1\right)  -\left(  H2\right)  $ hold and either $\left(  H3\right)  $ or
$\left(  H4\right)  $ is satisfied, then it is known \cite[Theorems 9 and
10]{KapKasVal14} that $\mathcal{A}_{\infty}$ is characterized by the union of
all bounded complete trajectories, that is,%
\begin{equation}
\mathcal{A}_{\infty}=\{\varphi\left(  0\right)  :\varphi\text{ is a bounded
complete trajectory of }\mathcal{R}_{0}\}. \label{CharactAttrM}%
\end{equation}

We will say that the family $\mathcal{R}$ satisfying $\left(  K1\right)
-\left(  K2\right)  $ is asymptotically autonomous if there is a family
$\mathcal{R}_{0}$ satisfying $\left(  H1\right)  -\left(  H2\right)  $ such
that for any sequences $\tau_{n}\rightarrow\infty,\ \varphi_{n}\in
\mathcal{R}_{\tau_{n}}$ such that $\varphi_{n}(\tau_{n})\rightarrow\varphi
_{0}$ there exist a subsequence $\{\psi_{n_{k}}\}$ of $\psi_{n}($%
\textperiodcentered$)=\varphi_{n}(\tau_{n}+$\textperiodcentered$)$ and
$\psi_{0}\in\mathcal{R}_{0}$ such that
\[
\psi_{n_{k}}(t)\rightarrow\psi_{0}(t)\text{ for all }t\geq0.
\]

\begin{lemma}
\label{ConvergAutonM}Assume that $\mathcal{R}=\{\mathcal{R}_{\tau}\}_{\tau
\in\mathbb{R}}$ is a family satisfying $\left(  K1\right)  -\left(  K2\right)
$ and such that the associated process $U$ has the pullback attractor
$\mathcal{A}=\{\mathcal{A}\left(  t\right)  \}$ and is forward asymptotically
compact. Also, let $\mathcal{R}$ be asymptotically autonomous and let the
multivalued semiflow corresponding to the limit family $\mathcal{R}_{0}$ have
the global attractor $\mathcal{A}_{\infty}$. Then%
\[
\lim_{t\rightarrow+\infty}dist\left(  \mathcal{A}\left(  t\right)
,\mathcal{A}_{\infty}\right)  =0.
\]

\end{lemma}

\begin{proof}
Let the statement be false. Then there are $\delta>0$, $t_{n}\nearrow+\infty$
and $y_{n}\in\mathcal{A}(t_{n})$ such that%
\[
dist(y_{n},\mathcal{A}_{\infty})\geq\delta\text{ for all }n.
\]
The inclusion $\mathcal{A}(t)\subset U(t,t_{0},\mathcal{A}(t_{0}))$ implies
that $\mathcal{A}$ is asymptotically compact. By Lemma \ref{ForwardA} we have
that $dist(\mathcal{A}(t),\omega(\mathcal{A}))\rightarrow0$ and then
$C_{0}=\overline{\cup_{t\geq\tau}\mathcal{A}(t)}$ is bounded for some $\tau$.
Since $\mathcal{A}_{\infty}$ attracts $C_{0}$, there is $N\left(
\delta\right)  $ for which%
\[
dist(G(t_{n},C_{0}),\mathcal{A}_{\infty})\leq\frac{\delta}{3}\text{ if }n\geq
N.
\]
As the pullback attractor is negatively invariant, there exist $z_{n}%
\in\mathcal{A}(t_{n}-t_{N})$ such that $y_{n}\in U(t_{n},t_{n}-t_{N},z_{n})$.
Passing to a subsequence we have that $z_{n}\rightarrow z_{0}\in C_{0}$. Let
$\varphi_{n}\in\mathcal{R}_{\tau_{n}}$, where $\tau_{n}=t_{n}-t_{N}$, be such
that $\varphi_{n}\left(  \tau_{n}\right)  =z_{n}$ and $\varphi_{n}%
(t_{n})=y_{n}$. Since $\mathcal{R}$ is asymptotically autonomous, there exists
$\varphi_{0}\in\mathcal{R}_{0}$ satisfying $\varphi_{0}(0)=z_{0}$ and a
subsequence of $\{\varphi_{n}\}$ (denoted the same) such that $\varphi
_{n}(t+\tau_{n})\rightarrow\varphi_{0}(t)$ for any $t\geq0$. Hence,%
\[
y_{n}=\varphi_{n}(\tau_{n}+t_{N})\rightarrow\varphi_{0}(t_{N}).
\]
Therefore, there is $N_{1}\geq N$ such that
\[
\rho(y_{n},\varphi_{0}(t_{N}))\leq\frac{\delta}{3}\text{ if }n\geq N_{1}.
\]
As $\varphi_{0}(t_{N})\in G(t_{N},z_{0})\subset$ $G(t_{N},C_{0})$, we derive
that%
\[
dist\left(  y_{n},\mathcal{A}_{\infty}\right)  \leq\rho(y_{n},\varphi
_{0}(t_{N}))+dist(G(t_{N},C_{0}),\mathcal{A}_{\infty})\leq\frac{2\delta}{3},
\]
which is a contradiction.
\end{proof}

\begin{remark}
This lemma is a version of Theorem 10 in \cite{LSSV} with slightly different conditions.
\end{remark}

\bigskip

\begin{theorem}
\label{ForwardAttractorTh2M}~Assume that $\mathcal{R}=\{\mathcal{R}_{\tau
}\}_{\tau\in\mathbb{R}}$ is a family satisfying $\left(  K1\right)  -\left(
K2\right)  $ and such that the associated process $U$ has the pullback
attractor $\mathcal{A}=\{\mathcal{A}\left(  t\right)  \}$ and is forward
asymptotically compact. Also, let $\mathcal{R}$ be asymptotically autonomous
and let the multivalued semiflow corresponding to the limit family
$\mathcal{R}_{0}$ have the global attractor $\mathcal{A}_{\infty}$. Assume
that%
\begin{equation}
dist_{H}(\mathcal{A}(t),\mathcal{A}_{\infty})\rightarrow0\text{ as
}t\rightarrow+\infty. \label{CondAInfM}%
\end{equation}
Then $\mathcal{A}$ is a forward attractor.
\end{theorem}

\begin{proof}
If we prove (\ref{CondOmegaM}) and (\ref{CondOmega2M}), then the result is a
consequence of Theorem \ref{ForwardAttractorThM}.

For assumption (\ref{CondOmega2M}) we only need to prove that $dist(\omega
(\mathcal{A}),\mathcal{A}(t))\rightarrow0$, because the other convergence
follows from Lemma \ref{ForwardA}. We will check that $\omega(\mathcal{A}%
)\subset\mathcal{A}_{\infty}$. Indeed, if $y\in\omega(\mathcal{A})$, there are
$y_{n}\in\mathcal{A}(t_{n})$, $t_{n}\rightarrow+\infty$, such that
$y_{n}\rightarrow y$. From (\ref{CondAInfM}) we have $dist(y_{n}%
,\mathcal{A}_{\infty})\rightarrow0$ and, therefore, the compactness of
$\mathcal{A}_{\infty}$ implies that $y\in\mathcal{A}_{\infty}$. Then, using
again (\ref{CondAInfM}) we conclude that%
\[
dist(\omega(\mathcal{A}),\mathcal{A}(t))\leq dist(\mathcal{A}_{\infty
},\mathcal{A}(t))\rightarrow0.
\]

Let us prove (\ref{CondOmegaM}). From (\ref{CondOmega2M}) and (\ref{CondAInfM}%
) we see that $\omega(\mathcal{A})=\mathcal{A}_{\infty}$. For $y_{0}\in
\omega_{f}(B,t_{0})$, there exist $y_{n}=u_{n}(t_{n})$, where $u_{n}%
\in\mathcal{R}_{t_{0}}$, with $t_{n}\rightarrow+\infty,\ x_{n}\in B$, such
that $y_{n}\rightarrow y_{0}$. By property (K2) we have that $u_{n}%
\in\mathcal{R}_{t_{0}+t_{n}}$. Let $v_{n}(t)=u_{n}(t+t_{n})$. Then
$v_{n}(0)\rightarrow y_{0}$ and since $\mathcal{R}$ is asymptotically
autonomous, there exists $v_{0}\in\mathcal{R}_{0}$ such that, up to a
subsequence, we have
\[
v_{n}(t)\rightarrow v_{0}(t)\text{ for any }t\geq0.
\]
As $v_{n}(t)=u_{n}(t+t_{n})\in U(t+t_{n},t_{0},x_{n}),$ we obtain that
$v_{0}(t)\in\omega_{f}(B,t_{0}).$ Put $\psi_{0}(t)=v_{0}(t)$ for $t\geq0$.
Further, let $v_{n}^{-1}(t)=u_{n}(t+t_{n}-1)\in U(t_{n}+t-1,t_{0},x_{n})$.
Then, since $U$ is forward asymptotically compact, up to a subsequence
$v_{n}^{-1}(0)\rightarrow y^{-1}\in\omega_{f}(B,t_{0})$. By property (K2) we
have that $u_{n}\in\mathcal{R}_{t_{0}+t_{n}-1}$. As $\mathcal{R}$ is
asymptotically autonomous, there exists $v_{-1}\in\mathcal{R}_{0}$ such that
\[
v_{n}^{-1}(t)\rightarrow v_{-1}(t)\text{ for any }t\geq0,
\]
As before $v_{-1}(t)\in\omega_{f}(B,t_{0}).$ We put $\psi_{-1}(t)=v_{-1}%
(t+1)$. Then $\psi_{-1}(t)=\psi_{0}(t),$ for any $t\geq0,$ and by (H2) we have
that $\psi_{-1}($\textperiodcentered$+s)\in\mathcal{R}_{0}$ for any $s\geq-1.$

Proceeding in this way for $k=-2,-3,...$ we obtain a sequence of functions
$\psi_{-k}:[-k,\infty)\rightarrow X$, $k=0,1,2,...,$ such that $\psi
_{-k}(t)\in\omega_{f}(B,t_{0})$, for $t\geq-k$, $\psi_{-k}(t)=\psi_{-k+1}(t),$
for any $t\geq-k+1,$ and $\psi_{-k}($\textperiodcentered$+s)\in\mathcal{R}%
_{0}$ for any $s\geq-k$.

We define $\psi:\left(  -\infty,+\infty\right)  \rightarrow X$ by $\psi\left(
t\right)  =\psi_{-k}(t)$, for $t\geq-k$, for any $k=0,1,2,...$ Then $\psi$ is
a bounded complete trajectory of $\mathcal{R}_{0}$ and $\psi\left(  0\right)
=y_{0}$. The characterization (\ref{CharactAttrM}) of the global attractor
implies then that $y_{0}\in\mathcal{A}_{\infty}=\omega(\mathcal{A})$. Thus,
$\omega_{f}(B,t_{0})\subset\omega(\mathcal{A})$.
\end{proof}

\bigskip

\section{Applications}

In this section, we will apply the abstract theory in order to prove that the
pullback attractor is a forward one for the Chafee-Infante equation, an
ordinary differential inclusion and a reaction-diffusion equation with
discontinuous nonlinearity. When the equations are asymptotically stable we
establish that the pullback attractor is continuous with respect to the
autonomous limit attractor as $t\rightarrow+\infty.$ Although this result is
used to prove that the existence of the forward attractor, it is interesting
by itself.

Throughout this section we will use the following notation. Let $H=L^{2}(0,1)$
and $V=H_{0}^{1}\left(  0,1\right)  $ with norms $\left\Vert
\text{\textperiodcentered}\right\Vert $ and $\left\Vert
\text{\textperiodcentered}\right\Vert _{V}$, respectively. An element $v\in H$
is said to be non-negative (denoted by $v\geq0$) if $v\left(  x\right)  \geq0$
for a.a. $x\in\left(  0,1\right)  $. An element $v\in V$ is said to be
positive (denoted by $v>0$) if $v(x)>0$ for all $x\in\left(  0,1\right)  $.

\subsection{The Chafee-Infante equation}

We consider the nonautonomous Chafee-Infante problem
\begin{equation}
\left\{
\begin{array}
[c]{l}%
\dfrac{\partial u}{\partial t}-\dfrac{\partial^{2}u}{\partial x^{2}}=\lambda
u-b(t)u^{3},\text{ on }(\tau,\infty)\times(0,\pi),\\
u(t,0)=u(t,1)=0,\\
u(\tau,x)=u_{\tau}(x),\text{ }x\in\left(  0,\pi\right)  ,
\end{array}
\right.  \label{Chafee}%
\end{equation}
where
\begin{equation}
1<\lambda<4 \label{CondLambda}%
\end{equation}
and $b:\mathbb{R}\rightarrow\mathbb{R}^{+}$ is an uniformly continuous and
differentiable function such that
\begin{equation}
0<b_{0}\leq b\left(  t\right)  \leq b_{1}. \label{Condb}%
\end{equation}

We will state first several results about the structure of the pullback
attractor from \cite{BroCarVal} and \cite{CLRBook}.

We recall that the eigenvalues of the operator $A=-\frac{\partial^{2}%
}{\partial x^{2}}$ with Dirichlet boundary condition on $\left(  0,\pi\right)
$ are $\lambda_{n}=n^{2}$, $n\geq1.$ We denote by $V^{2r}$ the spaces
$V^{2r}=D\left(  A^{r}\right)  $ for $r\in\mathbb{R}.$

For any $u_{\tau}\in V$ and $\tau\in\mathbb{R}$ there is a unique mild
solution $u\in C([\tau,\infty),V)$ to problem (\ref{Chafee}). The map
$U:\mathbb{R}_{\geq}^{2}\times V\rightarrow V$ given by $U(t,\tau,u_{\tau
})=u\left(  t\right)  $ is a continuous process. This procees prosseses a
pulback attractor $\mathcal{A}=\{\mathcal{A}(t)\}$ which satisfies that
$\cup_{t\in\mathbb{R}}\mathcal{A}(t)$ is bounded in $V$ and%
\[
\mathcal{A}(t)=\{\xi\left(  t\right)  :\xi\text{ is a bounded complete
trajectory}\}\text{.}%
\]
Also, there exists a maximal bounded complete trajectory $\xi_{M}^{+}$ such
that%
\[
-\xi_{M}^{+}\left(  t\right)  \leq\xi\left(  t\right)  \leq\xi_{M}^{+}\left(
t\right)  \text{ for any }t\in\mathbb{R},
\]
where $\xi$ is an arbitrary bounded complete trajectory \cite[Theorem
13.8]{CLRBook}. This obviously implies that $-\xi_{M}^{+}\left(  t\right)
\leq y\leq\xi_{M}^{+}\left(  t\right)  $ for any $y\in\mathcal{A}(t).$ The
functions $\xi_{M}^{+},\xi_{M}^{-}=-\xi_{M}^{+}$ are said to be nonautonomous
equilibria. Moreover, denote by $v_{1,b_{0}}^{+},v_{1,b_{1}}^{+}$ the positive
equilibria of the autonomous problem (\ref{Chafee}) with $b\left(  t\right)
\equiv b_{0}$ and $b\left(  t\right)  \equiv b_{1}$, respectively. Then%
\begin{equation}
v_{1,b_{0}}^{+}\leq\xi_{M}^{+}\left(  t\right)  \leq v_{1,b_{1}}^{+}\text{ for
any }t\in\mathbb{R}. \label{BoundsEquilibrium}%
\end{equation}
This implies the existence of $\varphi\in V$ such that $\varphi\left(
x\right)  >0$, for $x\in\left(  0,1\right)  $, and $\varphi\left(  x\right)
\leq\xi_{M}^{+}\left(  t,x\right)  $ for all $x\in\left(  0,1\right)  $,
$t\in\mathbb{R}$. In particular, $\xi_{M}^{+}$ is non-degenerate at
$t\rightarrow\pm\infty$. It is the unique bounded complete trajectory that is
non-degenerate at $-\infty$. The structure of the global attractor was
described in \cite[Section 4]{BroCarVal}, showing that any bounded complete
trajectory $\xi$ distinct from $0$ and $\xi_{M}^{\pm}$ has to be exclusively
of one of the following two types:

\begin{itemize}
\item $\xi\left(  t\right)  >0,$ for all $t\in\mathbb{R}$, and the following
convergences hold:%
\[
\xi\left(  t\right)  \rightarrow0\text{ as }t\rightarrow-\infty,
\]%
\[
\left\Vert \xi\left(  t\right)  -\xi_{M}^{+}\left(  t\right)  \right\Vert
_{V}\rightarrow0\text{ as }t\rightarrow+\infty.
\]

\item $\xi\left(  t\right)  <0,$ for all $t\in\mathbb{R}$, and the following
convergences hold:%
\[
\xi\left(  t\right)  \rightarrow0\text{ as }t\rightarrow-\infty,
\]%
\[
\left\Vert \xi\left(  t\right)  -\xi_{M}^{-}\left(  t\right)  \right\Vert
_{V}\rightarrow0\text{ as }t\rightarrow+\infty.
\]

\end{itemize}

Thus, the pullback attractor consists of the nonautonomous equilibria
$0,\xi_{M}^{+},\xi_{M}^{-}$ and the heteroclinic connections between them.
Finally, we recall that each set $\mathcal{A}(t)$ is connected \cite[Corollary
2.5]{CLRBook}.

We denote $V^{+}=\{v\in V:v\geq0\},\ V^{-}=\{v\in V:v\leq0\}$ and
$\mathcal{A}^{\pm}(t)=\mathcal{A}(t)\cap V^{\pm}.$

\begin{lemma}
\label{ConvergSequence}For any $t\in\mathbb{R}$ there exists a sequence
$\{y_{n}\}\subset\mathcal{A}^{+}(t)$ ($\subset\mathcal{A}^{-}(t)$),
$y_{n}\not =0$, such that $y_{n}\rightarrow0$.
\end{lemma}

\begin{proof}
Let us prove the result for $\mathcal{A}^{+}(t)$. It is clear that $0\leq
y\leq\xi_{M}^{+}\left(  t\right)  $ for any $t\in\mathbb{R}$, $y\in
\mathcal{A}^{+}(t)$. By contradiction assume that the statement is false. Then
there exists $\varepsilon>0$ such that any $z\not =0$ satisfying
$z\geq0,\ \left\Vert z\right\Vert <\varepsilon$ is not in $\mathcal{A}^{+}%
(t)$. We define the sets
\[
K_{1}=\mathcal{A}^{-}(t),\ K_{2}=\{z\in\mathcal{A}^{+}(t):\left\Vert
z\right\Vert \geq\varepsilon\}.
\]
Under our assumption it is obvious that $\mathcal{A}(t)=K_{1}\cup K_{2}$.
Also, there are disjoint open sets $U_{1},U_{2}$ such that $K_{i}\subset
U_{i}$, $i=1,2$. This implies that $\mathcal{A}(t)$ is not connected, which is
a contradiction.
\end{proof}

\begin{lemma}
\label{ForwardACChafee}$U$ is forward asymptotically compact.
\end{lemma}

\begin{proof}
Let $u$ be an arbitrary solution to problem (\ref{Chafee}) with $u_{\tau}\in
B$, a bounded set of $V$. Let $\alpha>\lambda-\lambda_{1}$ and $\gamma
=\alpha+\lambda_{1}-\lambda$. By $\alpha u^{2}\leq b_{0}u^{4}+\frac{\alpha
^{2}}{4b_{0}}$ we obtain in a standard way that%
\[
\frac{d}{dt}\left\Vert u\right\Vert ^{2}+2\gamma\left\Vert u\right\Vert
^{2}\leq\frac{d}{dt}\left\Vert u\right\Vert ^{2}+2\left(  \alpha
-\lambda\right)  \left\Vert u\right\Vert ^{2}+2\left\Vert u\right\Vert
_{V}^{2}\leq\frac{\alpha^{2}}{2b_{0}},
\]%
\begin{equation}
\left\Vert u(t)\right\Vert ^{2}\leq e^{-2\gamma(t-\tau)}\left\Vert u_{\tau
}\right\Vert ^{2}+\frac{\alpha^{2}}{4\gamma b_{0}}\ \forall t\geq\tau,
\label{Ineq1Chafee}%
\end{equation}%
\begin{equation}
\int_{t}^{t+r}\left\Vert u\right\Vert _{V}^{2}ds\leq\frac{1}{2}\left\Vert
u(t)\right\Vert ^{2}+\frac{\alpha^{2}}{4b_{0}}r\leq\frac{1}{2}e^{-2\gamma
(t-\tau)}\left\Vert u_{\tau}\right\Vert ^{2}+\frac{\alpha^{2}}{4\gamma b_{0}%
}\left(  \frac{1}{2}+r\right)  \text{ }\forall t\geq\tau,\ r>0.
\label{Ineq2Chafee}%
\end{equation}
Hence, there is $T\left(  B\right)  $ such that $\left\Vert u(t)\right\Vert
\leq\sqrt{1+\frac{\alpha^{2}}{4\gamma b_{0}}}=R_{0}$ for $t\geq T$. We
multiply the equation by $-\frac{\partial^{2}u}{\partial x^{2}}$. Then%
\begin{equation}
\frac{d}{dt}\left\Vert u\right\Vert _{V}^{2}\leq2\left(  \lambda-\lambda
_{1}\right)  \left\Vert u\right\Vert _{V}^{2}-6b\left(  t\right)  \int_{0}%
^{1}u^{2}u_{x}^{2}dx\leq2\left(  \lambda-\lambda_{1}\right)  \left\Vert
u\right\Vert _{V}^{2}. \label{Ineq3Chafee}%
\end{equation}
Integrating over $\left(  s,t+1\right)  $ we have%
\[
\left\Vert u(t+1)\right\Vert _{V}^{2}\leq\left\Vert u(s)\right\Vert _{V}%
^{2}+2\left(  \lambda-\lambda_{1}\right)  \int_{s}^{t+1}\left\Vert
u(r)\right\Vert _{V}^{2}dr.
\]
Integrating over $\left(  t,t+1\right)  $ and using (\ref{Ineq2Chafee}) we
obtain%
\begin{align}
\left\Vert u\left(  t+1\right)  \right\Vert _{V}^{2}  &  \leq\int_{t}%
^{t+1}\left\Vert u(r)\right\Vert _{V}^{2}dr+2\left(  \lambda-\lambda
_{1}\right)  \int_{t}^{t+1}\left\Vert u(r)\right\Vert _{V}^{2}dr\nonumber\\
&  \leq(1+2\left(  \lambda-\lambda_{1}\right)  )\left(  \frac{1}{2}\left\Vert
u(t)\right\Vert ^{2}+\frac{\alpha^{2}}{4b_{0}}\right) \nonumber\\
&  \leq(1+2\left(  \lambda-\lambda_{1}\right)  )\left(  \frac{R_{0}^{2}}%
{2}+\frac{\alpha^{2}}{4b_{0}}\right)  =R_{1}^{2}\text{ if }t\geq T.
\label{AcotNormV}%
\end{align}
We define the operator $F:\mathbb{R}\times V\rightarrow H$ given by
$F(u)=-b\left(  t\right)  u^{3}$. We see by $V\subset L^{6}(0,1)$ that%
\begin{equation}
\left\Vert F\left(  u\right)  \right\Vert =b\left(  t\right)  \left\Vert
u\right\Vert _{L^{6}}^{3}\leq C_{1}\left\Vert u\right\Vert _{V}^{3}.
\label{AcotF}%
\end{equation}
Then the variation of constants formula and (\ref{AcotNormV}) gives for $r<1$
that%
\begin{align*}
\left\Vert u\left(  t+2\right)  \right\Vert _{V^{2r}}  &  \leq\left\Vert
A^{r}e^{-A}u(t+1)\right\Vert +\int_{t+1}^{t+2}\left\Vert A^{r}e^{-A(t+2-s)}%
F(u(s))\right\Vert ds\\
&  \leq C_{2}+C_{3}\int_{t+1}^{t+2}(t+2-s)^{-r}\left\Vert u(s)\right\Vert
_{V}^{3}ds\\
&  \leq C_{2}+C_{3}R_{1}^{3}\frac{1}{1-r}\text{ if }t\geq T,
\end{align*}
where we have used the well-known inequality $\left\Vert A^{r}e^{-At}%
z\right\Vert \leq M_{r}t^{-r}e^{-at}\left\Vert z\right\Vert $ for some
constants $M_{r},a>0$ \cite{SellBook}$.$

From this we obtain that any sequence $y_{n}\in U(t_{n},\tau,B)$, where
$t_{n}\rightarrow+\infty$, is bounded in the space $V^{2r}$, which is
compactly embedded in $V$ for $r>\frac{1}{2}$. Therefore, $\{y_{n}\}$ is
relatively compact in $V$, proving the assertion.
\end{proof}

\bigskip

Let un consider the situation when our problem is asymptotically autonomous,
that is, there is $b_{0}\leq b\leq b_{1}$ such that%
\begin{equation}
b\left(  t\right)  \rightarrow b,\text{ as }t\rightarrow+\infty.
\label{CondAAutChafee}%
\end{equation}
It is well known \cite{Henry} that the autonomous limit problem, that is, the
one with $b\left(  t\right)  \equiv b$, generates a continuous semigroup
$S:\mathbb{R}^{+}\times V\rightarrow V$ having the global attractor
$\mathcal{A}_{\infty}$. Under assumption (\ref{CondLambda}) there are three
fixed points: $0,v_{1}^{+},v_{1}^{-}$, where $v_{1}^{+}>0$ and $v_{1}%
^{-}=-v_{1}^{+}$. The attractor consists of these fixed points and two bounded
complete trajectories $\varphi_{0}^{+},\varphi_{0}^{-}$ satisfying:

\begin{itemize}
\item $\varphi_{0}^{+}\left(  t\right)  >0,$ for all $t\in\mathbb{R}$, and the
following convergences hold:%
\[
\varphi_{0}^{+}\left(  t\right)  \rightarrow0\text{ as }t\rightarrow-\infty,
\]%
\[
\varphi_{0}^{+}\left(  t\right)  \rightarrow v_{1}^{+}\text{ as }%
t\rightarrow+\infty.
\]

\item $\varphi_{0}^{-}\left(  t\right)  <0,$ for all $t\in\mathbb{R}$, and the
following convergences hold:%
\[
\varphi_{0}^{-}\left(  t\right)  \rightarrow0\text{ as }t\rightarrow-\infty,
\]%
\[
\varphi_{0}^{-}\left(  t\right)  \rightarrow v_{1}^{-}\text{ as }%
t\rightarrow+\infty.
\]

\end{itemize}

\begin{lemma}
\label{AAChafee}The process $U$ is asymptotically autonomous.
\end{lemma}

\begin{proof}
Take $u_{\tau_{n}}\rightarrow u_{0}$ and the solutions $u_{n}\left(  t\right)
=U(t+\tau_{n},\tau_{n},u_{\tau_{n}})$ where $\tau_{n}\rightarrow+\infty$. We
define $v_{n}(t)=u_{n}(t+\tau_{n})$, which are solutions to problem
(\ref{Chafee}) with $v_{n}(0)=u_{\tau_{n}}$, $b\left(
\text{\textperiodcentered}\right)  =b_{n}\left(  \text{\textperiodcentered
}\right)  =b\left(  \text{\textperiodcentered}+\tau_{n}\right)  $. Also, $u$
is the solution to the autonomous problem (\ref{Chafee}) with $u(0)=u_{0}$,
$b\left(  \text{\textperiodcentered}\right)  \equiv b$. The difference
$w_{n}=v_{n}-u$ satisfies%
\[
\frac{\partial w_{n}}{\partial t}-\dfrac{\partial^{2}w_{n}}{\partial x^{2}%
}=\lambda w_{n}-b_{n}\left(  s\right)  v_{n}^{3}+bu^{3}=\lambda w_{n}%
+(b-b_{n}\left(  s\right)  )u^{3}+b_{n}\left(  s\right)  \left(  u^{3}%
-v_{n}^{3}\right)  .
\]
We multiply by $-\dfrac{\partial^{2}w_{n}}{\partial x^{2}}$ and use
\[
\int_{0}^{1}u^{3}(-\frac{\partial^{2}w_{n}}{\partial x^{2}})dx=3\int_{0}%
^{1}u^{2}\frac{\partial u}{\partial x}\frac{\partial w_{n}}{\partial x}%
dx\leq3\left\Vert u\right\Vert _{L^{\infty}}^{2}\left\Vert u\right\Vert
_{V}\left\Vert w_{n}\right\Vert _{V},
\]%
\begin{align*}
b_{n}(s)\int_{0}^{1}\left(  u^{3}-v_{n}^{3}\right)  \left(  -\frac
{\partial^{2}w_{n}}{\partial x^{2}}\right)  dx  &  =3b_{n}(s)\int_{0}%
^{1}\left(  u^{2}\frac{\partial u}{\partial x}-v_{n}^{2}\frac{\partial v_{n}%
}{\partial x}\right)  \frac{\partial w_{n}}{\partial x}dx\\
&  =-3b_{n}(s)\int_{0}^{1}\left(  u^{2}\left(  \frac{\partial w_{n}}{\partial
x}\right)  ^{2}+\left(  u+v_{n}\right)  w_{n}\frac{\partial v_{n}}{\partial
x}\frac{\partial w_{n}}{\partial x}\right)  dx\\
&  \leq3b_{n}\left(  s\right)  \left(  \left\Vert u\right\Vert _{L^{\infty}%
}^{2}\left\Vert w_{n}\right\Vert _{V}^{2}+\left(  \left\Vert u\right\Vert
_{L^{\infty}}+\left\Vert v_{n}\right\Vert _{L^{\infty}}\right)  \left\Vert
w_{n}\right\Vert _{L^{\infty}}\left\Vert v_{n}\right\Vert _{V}\left\Vert
w_{n}\right\Vert _{V}\right) \\
&  \leq Cb_{1}\left(  \left\Vert u\right\Vert _{L^{\infty}}^{2}+\left(
\left\Vert u\right\Vert _{L^{\infty}}+\left\Vert v_{n}\right\Vert _{L^{\infty
}}\right)  \left\Vert v_{n}\right\Vert _{V}\right)  \left\Vert w_{n}%
\right\Vert _{V}^{2}%
\end{align*}
in order to derive that%
\[
\frac{1}{2}\frac{d}{dt}\left\Vert w_{n}\right\Vert _{V}^{2}\leq\left(
\lambda-\lambda_{1}\right)  \left\Vert w_{n}\right\Vert _{V}^{2}+3\left\vert
b-b_{n}\left(  s\right)  \right\vert \left\Vert u\right\Vert _{L^{\infty}}%
^{2}\left\Vert u\right\Vert _{V}\left\Vert w\right\Vert _{V}+Cb_{1}\left(
\left\Vert u\right\Vert _{L^{\infty}}^{2}+\left(  \left\Vert u\right\Vert
_{L^{\infty}}+\left\Vert v_{n}\right\Vert _{L^{\infty}}\right)  \left\Vert
v_{n}\right\Vert _{V}\right)  \left\Vert w_{n}\right\Vert _{V}^{2}.
\]
The solutions $v_{n}(t),\ u(t)$ are uniformly bounded for $t\in\lbrack0,T]$,
$T>0$, in $V\subset L^{\infty}(0,1)$ (this follows from (\ref{Ineq3Chafee})
and Gronwall's lemma). Hence, there are constants $\alpha,\beta>0$ such that%
\[
\frac{d}{dt}\left\Vert w_{n}\right\Vert _{V}^{2}\leq\alpha\left\Vert
w_{n}\right\Vert _{V}^{2}+\beta\left\vert b-b_{n}\left(  s\right)  \right\vert
^{2}\text{ for }0<s<T.
\]
Therefore, for any $t\in\lbrack0,T]$ using (\ref{CondAAutChafee})\ we have%
\begin{align*}
\left\Vert w_{n}\left(  t\right)  \right\Vert _{V}^{2}  &  \leq e^{\alpha
t}\left\Vert u_{\tau}^{n}-u_{0}\right\Vert _{V}^{2}+\beta\int_{0}^{t}%
e^{\alpha\left(  t-s\right)  }\left\vert b-b_{n}\left(  s\right)  \right\vert
^{2}ds\\
&  \leq e^{\alpha t}\left\Vert u_{\tau}^{n}-u_{0}\right\Vert _{V}^{2}%
+\frac{\beta}{\alpha}e^{\alpha t}\sup_{s\geq0}\left\vert b-b_{n}\left(
s\right)  \right\vert ^{2}\rightarrow0\text{, as }n\rightarrow\infty.
\end{align*}
Since $T>0$ is arbitrary, the lemma is proved.
\end{proof}

\begin{lemma}
\label{ConvergEquilibriaChafee}$\xi_{M}^{\pm}(t)\rightarrow v_{1}^{\pm}$ as
$t\rightarrow+\infty.$
\end{lemma}

\begin{proof}
By Lemma \ref{ForwardACChafee}, $U$ is forward asymptotically compact, which
implies using Lemma \ref{ForwardOmega} that $\omega\left(  \xi_{M}%
^{+}(0)\right)  $ is non-empty, compact and attracts $\xi_{M}^{+}(0)$. Since
$dist(\xi_{M}^{+}\left(  t\right)  ,\omega\left(  \xi_{M}^{+}(0)\right)
)\rightarrow0$, it is enough to check that $\omega\left(  \xi_{M}%
^{+}(0)\right)  =v_{1}^{+}.$

Let $y_{0}\in\omega\left(  \xi_{M}^{+}(0)\right)  $. Then there is a sequence
$\{\tau_{n}\}$ such that $\xi_{M}^{+}\left(  \tau_{n}\right)  \rightarrow
y_{0}$, $\tau_{n}\rightarrow\infty.$ Let $u_{n}(t)=\xi_{M}^{+}\left(  \tau
_{n}+t\right)  =U(t+\tau_{n},\tau_{n},\xi_{M}^{+}\left(  \tau_{n}\right)  )$.
By Lemma \ref{AAChafee} we find that
\[
u_{n}\left(  t\right)  \rightarrow S\left(  t,y_{0}\right)  \text{ for any
}t\geq0.
\]
We set $\varphi_{0}\left(  t\right)  =S\left(  t,y_{0}\right)  $ for $t\geq0$.
Let now $u_{n}^{-1}\left(  t\right)  =\xi_{M}^{+}\left(  \tau_{n}+t-1\right)
=U(t+\tau_{n}-1,\tau_{n}-1,\xi_{M}^{+}\left(  \tau_{n}-1\right)  )$. Since by
Lemma \ref{ForwardACChafee} $U$ is forward asymptotically compact, up to a
subsequence $\xi_{M}^{+}\left(  \tau_{n}-1\right)  \rightarrow y^{-1}$. Again
using Lemma \ref{AAChafee} we deduce that
\[
u_{n}^{-1}(t)\rightarrow S(t,y^{-1})\text{ for any }t\geq0.
\]
We set $\varphi_{-1}\left(  t\right)  =S(t+1,y^{-1})$ for $t\geq-1$. It is
clear that $\varphi_{-1}\left(  t\right)  =\varphi_{0}\left(  t\right)  $ for
$t\geq0$. Also, $\varphi_{-1}(t)=S(t-s,\varphi_{-1}\left(  s\right)  )$ for
all $-1\leq s\leq t$. Proceeding in this same way for $k=2,3,...$ we obtain a
sequence $\varphi_{-k}$ such that
\begin{align*}
\varphi_{-k}\left(  t\right)   &  =\varphi_{-k+1}(t)\text{ for }t\geq-k+1,\\
\varphi_{-k}(t)  &  =S(t-s,\varphi_{-k}\left(  s\right)  )\text{ for }-k\leq
s\leq t.
\end{align*}
Let $\varphi$ be such that $\varphi\left(  t\right)  =\varphi_{-k}\left(
t\right)  $ for $t\geq-k$. The function $\varphi$ is a complete trajectory of
$S$. Moreover, it is bounded because $\varphi\left(  t\right)  \in
\omega\left(  \xi_{M}^{+}(0)\right)  ,$ for all $t\in\mathbb{R}$, and
(\ref{BoundsEquilibrium}) implies that $\varphi\left(  t\right)  \geq
v_{1,b_{0}}^{+}>0$ for any $t$. It follows that $\varphi\left(  t\right)
\equiv v_{1}^{+}$, proving that $\omega\left(  \xi_{M}^{+}(0)\right)
=v_{1}^{+}.$

The proof for $\xi_{M}^{-}$ is the same.
\end{proof}

\bigskip

Let us consider the operator $L=A-\lambda I$. The eigenvalues of $L$ are
$\overline{\lambda}_{n}=\lambda_{n}-\lambda$, $n\geq1$. From condition
(\ref{CondLambda}) we see that $\overline{\lambda}_{1}<0<\overline{\lambda
}_{2}<\overline{\lambda}_{3}<...$ Then $H=H_{1}\oplus H_{2}$, where $H_{1}$ is
the one-dimensional space generated by the eigenfunction for $\overline
{\lambda}_{1}$ and $H_{2}$ is the subspace generated by the eigenfunction for
$\{\overline{\lambda}_{2},\overline{\lambda}_{3},...\}$. Let $Q$ be the
projection onto $H_{1}$ and $P=I-Q$. The spaces $H_{j}$ are invariant for the
operator $L$. We denote $L_{j}=L\mid_{H_{j}}$. It is well known \cite[Theorem
1.5.3]{Henry} that%
\begin{equation}
\left\Vert e^{-L_{1}t}z\right\Vert _{V}\leq Ce^{-\overline{\lambda}_{1}%
t}\left\Vert z\right\Vert _{V}\text{ for }t\leq0, \label{IneqL1}%
\end{equation}%
\begin{equation}
\left\Vert e^{-L_{2}t}z\right\Vert _{V}\leq Ce^{-\overline{\lambda}_{2}%
t}\left\Vert z\right\Vert _{V}\text{ for }t\geq0. \label{IneqL2}%
\end{equation}

\begin{lemma}
The bounded complete trajectories $\varphi_{0}^{\pm}$ of the autonomous
problem satisfy%
\begin{equation}
\left\Vert \varphi_{0}^{\pm}\left(  \tau\right)  \right\Vert _{V}\leq
Ce^{-\overline{\lambda}_{1}\tau}\text{ }\forall\tau\leq\tau_{0},
\label{AcotFi0}%
\end{equation}
for some constants $C>0,\ \tau_{0}<0.$
\end{lemma}

\begin{proof}
We can write $\varphi_{0}^{+}$ as $\varphi_{0}^{+}\left(  t\right)  =p\left(
t\right)  +q\left(  t\right)  $, where $p\left(  t\right)  =P\varphi_{0}%
^{+}\left(  t\right)  ,\ q\left(  t\right)  =Q\varphi_{0}^{+}\left(  t\right)
$. If $F:V\rightarrow V$ is the operator defined by $F(u)=-bu^{3}$, then the
variation of constants formula gives%
\[
q\left(  \tau\right)  =e^{-L_{1}\tau}q\left(  0\right)  +\int_{0}^{\tau
}e^{-L_{1}(\tau-s)}QF\left(  \varphi_{0}^{+}\left(  s\right)  \right)
ds\text{ for }\tau\leq0.
\]
We observe that%
\[
\left\Vert F\left(  u\right)  \right\Vert _{V}^{2}=\int_{0}^{1}9u^{4}u_{x}%
^{2}dx\leq9\left\Vert u\right\Vert _{L^{\infty}}^{4}\left\Vert u\right\Vert
_{V}^{2}\leq K\left\Vert u\right\Vert _{V}^{6}.
\]
Then using (\ref{IneqL1}) we have%
\begin{align*}
\left\Vert q\left(  \tau\right)  \right\Vert _{V}  &  \leq C_{1}%
e^{-\overline{\lambda}_{1}\tau}+C_{2}\int_{0}^{\tau}e^{-\overline{\lambda}%
_{1}(\tau-s)}\left\Vert \varphi_{0}^{+}\left(  s\right)  \right\Vert _{V}%
^{3}ds\\
&  \leq C_{1}e^{-\overline{\lambda}_{1}\tau}+\frac{C_{2}}{\left\vert
\overline{\lambda}_{1}\right\vert }\sup_{s\in\mathbb{R}}\left\Vert \varphi
_{0}^{+}\left(  s\right)  \right\Vert _{V}^{3}e^{-\overline{\lambda}_{1}\tau
}\leq C_{3}e^{-\overline{\lambda}_{1}\tau}\text{ }\forall\tau\leq0.
\end{align*}

The map satisfies the properties:%
\[
F(0)=0,
\]%
\begin{align*}
\left\Vert F(u)-F(v)\right\Vert ^{2}  &  =b^{2}\int_{0}^{1}\left(  u^{3}%
-v^{3}\right)  ^{2}dx=b^{2}\int_{0}^{1}\left(  u-v\right)  ^{2}\left(
u^{2}+uv+v^{2}\right)  ^{2}dx\\
&  \leq b^{2}\left(  \left\Vert u\right\Vert _{L^{\infty}}^{2}+\left\Vert
u\right\Vert _{L^{\infty}}\left\Vert v\right\Vert _{L^{\infty}}+\left\Vert
v\right\Vert _{L^{\infty}}^{2}\right)  ^{2}\left\Vert u-v\right\Vert ^{2}\\
&  \leq C_{4}\left(  \left\Vert u\right\Vert _{V}^{2}+\left\Vert v\right\Vert
_{V}^{2}\right)  ^{2}\left\Vert u-v\right\Vert _{V}^{2}.
\end{align*}
Therefore, since $\varphi_{0}^{+}$ belongs to the unstable manifold of $0$,
the saddle-point property \cite[Theorem 5.2.1]{Henry} implies that%
\[
\left\Vert p\left(  \tau\right)  \right\Vert _{V}=o(\left\Vert q\left(
\tau\right)  \right\Vert _{V})\text{ as }\tau\rightarrow-\infty.
\]
Thus, there is $\tau_{0}<0$ such that
\[
\left\Vert p\left(  \tau\right)  \right\Vert _{V}\leq C_{3}e^{-\overline
{\lambda}_{1}\tau}\text{ }\forall\tau\leq\tau_{0}.
\]
The result follows.
\end{proof}

\bigskip

\begin{theorem}
\label{CondAInfChafee}$\lim_{t\rightarrow+\infty}dist_{H}\left(
\mathcal{A}\left(  t\right)  ,\mathcal{A}_{\infty}\right)  =0.$
\end{theorem}

\begin{proof}
Let us prove that $dist(\mathcal{A}_{\infty},\mathcal{A}\left(  t\right)
)\rightarrow0$ as $t\rightarrow+\infty.$

We will establish first this property with respect to the topology of the
space $H$. In this case, the Hausdorff semidistance from $A$ to $B$ will be
denoted by $dist_{L^{2}}\left(  A,B\right)  .$

Let $\varepsilon>0,\ z_{0}\in\mathcal{A}_{\infty}$ be arbitrary. First let
$z\geq0$. If $\left\Vert z_{0}\right\Vert \leq\varepsilon$, then $dist_{L^{2}%
}\left(  z_{0},\mathcal{A}\left(  t\right)  \right)  \leq\varepsilon$ for all
$t$. If $z_{0}\in\mathcal{A}_{\infty}$ is such that $\left\Vert z_{0}%
-v_{1}^{+}\right\Vert \leq\frac{\varepsilon}{2}$, then $\xi_{M}\left(
t\right)  \rightarrow v_{1}^{+}$ (see Lemma \ref{ConvergEquilibriaChafee})
implies that there is $T_{0}\left(  \varepsilon\right)  $ for which
\[
dist_{L^{2}}(z_{0},\mathcal{A}\left(  t\right)  )\leq\left\Vert z_{0}%
-v_{1}^{+}\right\Vert +\left\Vert v_{1}^{+}-\xi_{M}^{+}\left(  t\right)
\right\Vert \leq\varepsilon\text{ for any }t\geq T_{0}.
\]

Let now $\left\Vert z_{0}-z_{1}^{+}\right\Vert >\frac{\varepsilon}{2}$,
$\left\Vert z_{0}\right\Vert >\varepsilon$. The point $z_{0}$ belongs to the
bounded complete trajectory $\varphi_{0}^{+}$. By (\ref{AcotFi0}) we choose
$t_{\varepsilon},\ a_{0}(\varepsilon)$ such that
\[
\left\Vert \varphi_{0}^{+}\left(  t_{\varepsilon}\right)  \right\Vert \leq
Ce^{\left(  \lambda-\lambda_{1}\right)  t_{\varepsilon}}=\varepsilon,
\]%
\[
\varphi_{0}(t_{\varepsilon}+a_{0})=z_{0}.
\]
Since $\varphi_{0}\left(  t\right)  \overset{t\rightarrow+\infty}{\rightarrow
}v_{1}^{+}$, there is $a_{1}\left(  \varepsilon\right)  >$ such that
$\left\Vert \varphi_{0}^{+}\left(  t\right)  -v_{1}^{+}\right\Vert \leq
\frac{\varepsilon}{2}$ if $t\geq a_{1}+t_{\varepsilon}$. Hence, $a_{0}\leq
a_{1}$.

Put $\varphi_{\tau}(t)=\varphi_{0}^{+}(t+t_{\varepsilon}-a_{0}-\tau)$. Then
$\left\Vert \varphi_{\tau}(\tau)\right\Vert =\varphi_{0}^{+}(t_{\varepsilon
}-a_{0})$, $\varphi_{\tau}(\tau+2a_{0})=z_{0}$ and%
\begin{equation}
\left\Vert \varphi_{\tau}\left(  \tau\right)  \right\Vert \leq Ce^{\left(
\lambda-\lambda_{1}\right)  (t_{\varepsilon}-a_{0})}\leq\varepsilon
e^{-\left(  \lambda-\lambda_{1}\right)  a_{0}}. \label{IneqFitao}%
\end{equation}

We know that any $y\in\mathcal{A}\left(  t\right)  $ satisfies $0\leq y\leq
\xi_{M}(t)$. By Lemma \ref{ConvergSequence} for any $t_{0}$ there is
$y_{\varepsilon,t_{0}}\in\mathcal{A}\left(  t_{0}\right)  $ such that
\begin{equation}
\left\Vert y_{\varepsilon,t_{0}}\right\Vert \leq\varepsilon e^{-\left(
\lambda-\lambda_{1}\right)  a_{1}}. \label{IneqyyEps}%
\end{equation}
Take a bounded complete trajectory $\psi_{\varepsilon,t_{0}}$ such that
$\psi_{\varepsilon,t_{0}}\left(  t_{0}\right)  =y_{\varepsilon,t_{0}}$. For
$t$ arbitrarily big we set $t_{0}=t-2a_{0}$. The function $v\left(  s\right)
=\psi_{\varepsilon,t-2a_{0}}\left(  s\right)  $ is the solution to problem
(\ref{Chafee}) with $v\left(  t-2a_{0}\right)  =y_{\varepsilon,t_{0}}$. Also,
for $\tau=t-2a_{0}$, the function $u\left(  t\right)  =\varphi_{t-2a_{0}%
}\left(  t\right)  $ is the solution to problem (\ref{Chafee}) with $u\left(
t-2a_{0}\right)  =\varphi_{0}^{+}\left(  t_{\varepsilon}-a_{0}\right)  $ and
$b\left(  t\right)  \equiv b$. The difference $w=v-u$ satisfies%
\[
\frac{\partial w}{\partial t}-\frac{\partial^{2}w}{\partial x^{2}}=\lambda
w-b\left(  s\right)  v^{3}+bu^{3}=\lambda w+(b-b\left(  s\right)
)v^{3}+b\left(  s\right)  \left(  u^{3}-v^{3}\right)  .
\]
Multiplying by $w$ and using%
\[
\int_{0}^{1}v^{3}wdx\leq\left\Vert v\right\Vert _{L^{6}}^{3}\left\Vert
w\right\Vert ,
\]
we obtain%
\[
\frac{d}{dt}\left\Vert w\right\Vert ^{2}\leq(\lambda-\lambda_{1}+\frac
{1}{a_{1}(\varepsilon)})\left\Vert w\right\Vert ^{2}+a_{1}(\varepsilon
)\frac{\left\vert b-b\left(  s\right)  \right\vert ^{2}}{4}\left\Vert
v\right\Vert _{L^{6}}^{6}-b_{0}\int_{0}^{1}w^{2}\left(  u^{2}+uv+v^{2}\right)
dx.
\]

The functions$\ v\left(  t\right)  $ are uniformly bounded in $V\subset
L^{6}(0,1)$. Therefore, there is a constant $K>0$ such that
\[
\frac{d}{dt}\left\Vert w\right\Vert ^{2}\leq(\lambda-\lambda_{1}+\frac
{1}{a_{1}(\varepsilon)})\left\Vert w\right\Vert ^{2}+Ka_{1}\left(
\varepsilon\right)  \left\vert b-b\left(  s\right)  \right\vert ^{2}.
\]
We choose $T_{1}\left(  \varepsilon\right)  >0$ such that%
\[
\left\vert b\left(  s\right)  -b\right\vert \leq\frac{\varepsilon}{\sqrt
{a_{1}(\varepsilon)}}e^{-\left(  \lambda-\lambda_{1}\right)  a_{1}}\text{ for
any }s\geq T_{1}.
\]
Then, using (\ref{IneqFitao}), (\ref{IneqyyEps}) for $t\geq T_{1}+2a_{1}$ we
obtain%
\begin{align*}
\left\Vert z_{0}-\psi_{\varepsilon,t-2_{a_{0}}}(t)\right\Vert ^{2}  &
=\left\Vert w(t)\right\Vert ^{2}\leq e^{(\lambda-\lambda_{1}+\frac{1}{a_{1}%
})2a_{0}}\left\Vert w(t-2a_{0})\right\Vert ^{2}+\frac{Ka_{1}}{\lambda
-\lambda_{1}+\frac{1}{a_{1}}}e^{(\lambda-\lambda_{1}+\frac{1}{a_{1}})2a_{0}%
}\sup_{s\geq t-2a_{0}}\left\vert b-b\left(  s\right)  \right\vert ^{2}\\
&  \leq2e^{(\lambda-\lambda_{1}+\frac{1}{a_{1}})2a_{0}}\left(  \left\Vert
y_{\varepsilon,t-2_{a_{0}}}\right\Vert ^{2}+\left\Vert \varphi_{t-2a_{0}%
}\left(  t-2a_{0}\right)  \right\Vert ^{2}\right)  +\frac{Ka_{1}}%
{\lambda-\lambda_{1}}e^{(\lambda-\lambda_{1}+\frac{1}{a_{1}})2a_{0}}%
\sup_{s\geq t-2a_{1}}\left\vert b-b\left(  s\right)  \right\vert ^{2}\\
&  \leq e^{2}(4+\frac{K}{\lambda-\lambda_{1}})\varepsilon^{2}=R^{2}%
\varepsilon^{2}.
\end{align*}

Fo $z\leq0$ the proof is similar.

Hence, there is $T\left(  \varepsilon\right)  $ such that
\[
dist_{H}(\mathcal{A}_{\infty},\mathcal{A}\left(  t\right)  )=\sup_{z_{0}%
\in\mathcal{A}_{\infty}}dist_{L^{2}}(z_{0},\mathcal{A}\left(  t\right)  )\leq
R\varepsilon\text{ if }t\geq T.
\]
As $\varepsilon$ is arbitrary, we have proved that%
\begin{equation}
dist_{L^{2}}(\mathcal{A}_{\infty},\mathcal{A}\left(  t\right)  )\rightarrow
0\text{ as }t\rightarrow+\infty. \label{DistH}%
\end{equation}

Second, we will prove that $dist(\mathcal{A}_{\infty},\mathcal{A}\left(
t\right)  )\rightarrow0$ as $t\rightarrow+\infty.$

By contradiction, if this not true, then there are $\varepsilon>0$ and
sequences $y_{n}\in\mathcal{A}_{\infty}$, $t_{n}\rightarrow+\infty$ such that%
\[
dist\left(  y_{n},\mathcal{A}\left(  t_{n}\right)  \right)  \geq
\varepsilon\text{ }\forall n.
\]
Take $z_{n}\in\mathcal{A}\left(  t_{n}\right)  $ such that $dist_{L^{2}%
}\left(  y_{n},\mathcal{A}\left(  t_{n}\right)  \right)  =\left\Vert
y_{n}-z_{n}\right\Vert $. By the invariance of the pullback attractor we have
that $z_{n}\in U(t_{n},0,\mathcal{A}\left(  0\right)  )$. Since by Lemma
\ref{ForwardACChafee} $U$ is forward asymptotically compact, we obtain passing
to a subsequence that $z_{n}\rightarrow z_{0}$ in $V$ (and then in $H$ as
well). By the compactness of $\mathcal{A}_{\infty}$, we have $y_{n}\rightarrow
y_{0}\in\mathcal{A}_{\infty}$ in $V$. By (\ref{DistH}) we deduce that
$\left\Vert y_{n}-z_{n}\right\Vert \rightarrow0$, and then $y_{0}=z_{0}%
\in\mathcal{A}_{\infty}$. Thus,
\[
dist\left(  y_{n},\mathcal{A}\left(  t_{n}\right)  \right)  \leq\left\Vert
y_{n}-z_{n}\right\Vert _{V}\rightarrow0\text{ as }n\rightarrow\infty.
\]
This is a contradiction.

The fact that $\lim_{t\rightarrow+\infty}dist\left(  \mathcal{A}\left(
t\right)  ,\mathcal{A}_{\infty}\right)  =0$ follows from Lemmas \ref{AAChafee}%
, \ref{ForwardACChafee} and \ref{ConvergAuton}.
\end{proof}

\bigskip

\begin{theorem}
$\mathcal{A}$ is a forward attractor.
\end{theorem}

\begin{proof}
It follows from Lemmas \ref{AAChafee}, \ref{ForwardACChafee},
\ref{CondAInfChafee} and Theorem \ref{ForwardAttractorTh2}.
\end{proof}

\subsection{An ordinary differential inclusion}

We consider the problem%
\begin{equation}
\left\{
\begin{array}
[c]{c}%
\dfrac{du}{dt}+\lambda u\in b\left(  t\right)  H_{0}(u),\ t\geq\tau,\\
u(\tau)=u_{\tau},
\end{array}
\right.  \label{OrdIncl}%
\end{equation}
where $\lambda>0$, $b:\mathbb{R}\rightarrow\mathbb{R}^{+}$ is a continuous
functions satisfying%
\[
0<b_{0}\leq b\left(  t\right)  \leq b_{1}\text{, for }t\in\mathbb{R},
\]
and $H$ stands for the Heaviside function given by%
\begin{equation}
H_{0}(u)=\left\{
\begin{array}
[c]{c}%
1\text{ if }u>0,\\
\lbrack-1,1]\text{ if }u=0,\\
-1\text{ if }u<0.
\end{array}
\right.  \label{Heaviside}%
\end{equation}

We say that the function $u:[\tau,\infty)\rightarrow\mathbb{R}$ is a solution
to (\ref{OrdIncl}) if $u\in C([\tau,\infty),\mathbb{R}),\ \dfrac{du}{dt}\in
L_{loc}^{\infty}([s,\infty),\mathbb{R})$ and there exists $h\in L_{loc}%
^{\infty}([\tau,\infty),\mathbb{R})$ satisfying $h\left(  t\right)  \in
H_{0}(u(t))$, for a.a. $t\in\left(  \tau,\infty\right)  $, and%
\[
\dfrac{du}{dt}+\lambda u=b\left(  t\right)  h\left(  t\right)  ,\ \text{for
a.a.\ }t\geq\tau.
\]
In \cite[Corollary 3.5]{CLV16} it is shown that all the possible solutions to
problem (\ref{OrdIncl}) are the following:%
\begin{equation}
u\left(  t\right)  =e^{-\lambda\left(  t-\tau\right)  }u_{\tau}+\int_{\tau
}^{t}e^{-\lambda\left(  t-r\right)  }b\left(  r\right)  dr\text{ if }u_{\tau
}>0, \label{Sol1}%
\end{equation}%
\begin{equation}
u\left(  t\right)  =e^{-\lambda\left(  t-\tau\right)  }u_{\tau}-\int_{\tau
}^{t}e^{-\lambda\left(  t-r\right)  }b\left(  r\right)  dr\text{ if }u_{\tau
}<0, \label{Sol2}%
\end{equation}%
\begin{equation}
u_{\infty}\left(  t\right)  =0\text{ for all }t\geq\tau, \label{Sol3}%
\end{equation}%
\begin{equation}
u_{r}^{+}\left(  t\right)  =\left\{
\begin{array}
[c]{c}%
0\text{ si }\tau\leq t\leq r,\\
\int_{r}^{t}e^{-\lambda\left(  t-s\right)  }b\left(  s\right)  ds,\text{ if
}t\geq r,
\end{array}
\right.  \label{Sol4}%
\end{equation}%
\begin{equation}
u_{r}^{-}\left(  t\right)  =\left\{
\begin{array}
[c]{c}%
0\text{ si }\tau\leq t\leq r,\\
-\int_{r}^{t}e^{-\lambda\left(  t-s\right)  }b\left(  s\right)  ds,\text{ if
}t\geq r,
\end{array}
\right.  \label{Sol5}%
\end{equation}

The families $\mathcal{R}_{\tau}\subset C([\tau,\infty),\mathbb{R})$ will be
the set of all solutions $u$ to (\ref{OrdIncl}). It is shown in \cite{CLV16}
that $\mathcal{R}$ satisfies the properties (K1)-(K4) and that the
corresponding strict multivalued process $U$ possesses a strictly invariant
pullback attractor $\mathcal{A}=\{\mathcal{A}(t)\}$, which is globally
bounded, that is, $\cup_{t\in\mathbb{R}}\mathcal{A}(t)$ is bounded, and
satisfies%
\[
\mathcal{A}(t)=\{\gamma\left(  t\right)  :\gamma\text{ is a bounded complete
trajectory of }\mathcal{R}\}.
\]
Moreover, the structure of the pullback attractor was fully described in
\cite{CLV16}. Namely, it was shown that all the possible bounded complete
trajectories are the following ones:%
\[
\xi_{0}\left(  t\right)  \equiv0,
\]%
\begin{align*}
\xi_{M}^{+}(t)  &  =\int_{-\infty}^{t}e^{-\lambda\left(  t-s\right)  }b\left(
s\right)  ds,\\
\xi_{M}^{-}(t)  &  =-\int_{-\infty}^{t}e^{-\lambda\left(  t-s\right)
}b\left(  s\right)  ds,
\end{align*}%
\begin{align*}
u_{r}^{+}(t)  &  =\left\{
\begin{array}
[c]{c}%
0\text{ if }t\leq r,\\
\int_{r}^{t}e^{-\lambda\left(  t-s\right)  }b\left(  s\right)  ds\text{ if
}t\geq r,
\end{array}
\right. \\
u_{r}^{-}(t)  &  =\left\{
\begin{array}
[c]{c}%
0\text{ if }t\leq r,\\
-\int_{r}^{t}e^{-\lambda\left(  t-s\right)  }b\left(  s\right)  ds\text{ if
}t\geq r,
\end{array}
\right.
\end{align*}
with $r\in\mathbb{R}$ arbitrary. The function $\xi_{M}^{+}$ ($\xi_{M}^{-}$) is
the only bounded strictly positive (negative) complete trajectory. While
$\xi_{0}$ is a fixed point in the classical sense, $\xi_{M}^{+}$, $\xi_{M}%
^{-}$ are said to be nonautonomous stationary solutions. The functions
$\xi_{M}^{+}$, $\xi_{M}^{-}$ are upper and lower bounds of the attractor, that
is,%
\[
\xi_{M}^{-}\left(  t\right)  \leq\gamma\left(  t\right)  \leq\xi_{M}%
^{+}\left(  t\right)  \text{ for all }t\in\mathbb{R},
\]
where $\gamma$ is any bounded complete trajectory of $\mathcal{R}$. Moreover,
it is easy to see that%
\[
\mathcal{A}(t)=[\xi_{M}^{-}\left(  t\right)  ,\xi_{M}^{+}\left(  t\right)  ].
\]
The functions $u_{r}^{+},u_{r}^{+}$ connect these stationary solutions in the
sense that%
\begin{align*}
\left\vert u_{r}^{+}\left(  t\right)  -\xi_{M}^{+}(t)\right\vert  &
\rightarrow0\text{ as }t\rightarrow+\infty,\\
u_{r}^{+}\left(  t\right)   &  \rightarrow0\text{ as }t\rightarrow-\infty,
\end{align*}%
\begin{align*}
\left\vert u_{r}^{-}\left(  t\right)  -\xi_{M}^{-}(t)\right\vert  &
\rightarrow0\text{ as }t\rightarrow+\infty,\\
u_{r}^{-}\left(  t\right)   &  \rightarrow0\text{ as }t\rightarrow-\infty.
\end{align*}
Thus, the pullback attractor consists of the nonautonomous stationary
solutions $\xi_{0},\ \xi_{M}^{+}$, $\xi_{M}^{-}$ and their heteroclinic
connections $u_{r}^{+},\ u_{r}^{-}$.

\begin{lemma}
\label{ForwardACOrdIncl}$U$ is forward asymptotically compact.
\end{lemma}

\begin{proof}
It is straightforward to see from (\ref{Sol1})-(\ref{Sol5}) that $\cup
_{t\geq\tau}U(t,\tau,B)$ is a bounded set for any $B$ bounded and any $\tau
\in\mathbb{R}$. Therefore, any sequence $y_{n}\in U(t_{n},\tau,B)$,
$t_{n}\rightarrow+\infty$, is relatively compact.
\end{proof}

\bigskip

Let us prove that the pullback attractor is a forward attractor.

\begin{theorem}
$\mathcal{A}$ is a forward attractor.
\end{theorem}

\begin{proof}
Let $B$ be a bounded set and $R>0$ such that $\left\vert z\right\vert \leq R$
for any $z\in B$. For $\varepsilon>0,\ t_{0}\in\mathbb{R}$ take $T\left(
\varepsilon,t_{0}\right)  >t_{0}$ such that%
\[
e^{-\lambda\left(  t-t_{0}\right)  }R\leq\frac{\varepsilon}{2},\ e^{-\lambda
t}\int_{-\infty}^{t_{0}}e^{\lambda s}b\left(  s\right)  ds\leq\frac
{\varepsilon}{2}\text{ if }t\geq T.
\]
Let $y\in U(t,t_{0},x)\subset U(t,t_{0},B)$. If $x>0$, then by (\ref{Sol1}) we
have%
\begin{align*}
\left\vert y-\xi_{M}^{+}(t)\right\vert  &  =\left\vert e^{-\lambda\left(
t-t_{0}\right)  }x+\int_{t_{0}}^{t}e^{-\lambda\left(  t-s\right)  }b\left(
s\right)  ds-\int_{-\infty}^{t}e^{-\lambda\left(  t-s\right)  }b\left(
s\right)  ds\right\vert \\
&  \leq e^{-\lambda\left(  t-t_{0}\right)  }R+e^{-\lambda t}\int_{-\infty
}^{t_{0}}e^{\lambda s}b\left(  s\right)  ds\leq\varepsilon.
\end{align*}
If $x<0$,\ by (\ref{Sol2}) we obtain%
\begin{align*}
\left\vert y-\xi_{M}^{-}(t)\right\vert  &  =\left\vert e^{-\lambda\left(
t-t_{0}\right)  }x-\int_{t_{0}}^{t}e^{-\lambda\left(  t-s\right)  }b\left(
s\right)  ds+\int_{-\infty}^{t}e^{-\lambda\left(  t-s\right)  }b\left(
s\right)  ds\right\vert \\
&  \leq e^{-\lambda\left(  t-t_{0}\right)  }R+e^{-\lambda t}\int_{-\infty
}^{t_{0}}e^{\lambda s}b\left(  s\right)  ds\leq\varepsilon.
\end{align*}
If $x=0$, by by (\ref{Sol3})-(\ref{Sol5}) we have three possibilities:

\begin{enumerate}
\item $y=0;$

\item $y=u_{r}^{+}\left(  t\right)  $ for some $r\geq s;$

\item $y=u_{r}^{-}\left(  t\right)  $ for some $r\geq s.$
\end{enumerate}

Hence, $y\in\mathcal{A}\left(  t\right)  $. We deduce that%
\[
dist(U(t,t_{0},B),\mathcal{A}\left(  t\right)  )\leq\varepsilon\text{ if
}t\geq T.
\]
Thus, $\mathcal{A}$ is a forward attractor.
\end{proof}

\bigskip

We observe that condition (\ref{CondOmega0M}) is not satisfied in this
example. Indeed, take $b\left(  t\right)  =2+\sin(t)$, $\lambda=1$. Then
\[
\xi_{M}^{+}(t)=\int_{-\infty}^{t}e^{-\left(  t-s\right)  }\left(
2+\sin(s)\right)  ds=\frac{1}{2}\sin t-\frac{1}{2}\cos t+2
\]
and $\xi_{M}^{+}(t)$ oscilates in the interval $[2-a,2+a]$, where $a=\frac
{1}{2}\sin\frac{3\pi}{4}-\frac{1}{2}\cos\frac{3\pi}{4}=\frac{\sqrt{2}}{2}$.
Since $\mathcal{A}\left(  t\right)  =[\xi_{M}^{+}(t),-\xi_{M}^{+}(t)]$, we see
that%
\begin{align*}
\omega(\mathcal{A})  &  =[-2-a,2+a],\\
\omega_{0}(\mathcal{A})  &  =[-2+a,2-a].
\end{align*}
Hence, $\omega(\mathcal{A})\not =\omega_{0}(\mathcal{A})$ and then by Theorem
\ref{ForwardAttractorTh0M} condition (\ref{CondOmega0M}) is not true.

In the particular situation where problem (\ref{OrdIncl}) is asymptotically
autonomous, we can prove that (\ref{CondOmega0M}) is satisfied. Let
\[
b\left(  t\right)  \rightarrow b\in\lbrack b_{0},b_{1}]\text{ as }%
t\rightarrow+\infty.
\]
For the limit system with $b\left(  t\right)  \equiv b$, it is known
\cite{CLV16} that if $\mathcal{R}_{0}$ is the set of all solutions to
(\ref{OrdIncl}), then properties (H1)-(H4) are satisfied and the corresponding
strict multivalued semiflow $G$ has an invariant attractor $\mathcal{A}%
_{\infty}$. Let us establish condition (\ref{CondAInfM}).

In this case, it is shown in \cite{CLV16} that all the possible bounded
complete trajectories are given by the fixed points%
\[
z_{0}=0,\ z_{1}^{+}=\frac{b}{\lambda},\ z_{1}^{-}=-\frac{b}{\lambda},
\]
and the functions%
\begin{align*}
\varphi_{r}^{+}(t)  &  =\left\{
\begin{array}
[c]{c}%
0\text{ if }t\leq r,\\
\frac{b}{\lambda}(1-e^{-\lambda\left(  t-r\right)  })\text{ if }t\geq r,
\end{array}
\right. \\
\varphi_{r}^{-}(t)  &  =\left\{
\begin{array}
[c]{c}%
0\text{ if }t\leq r,\\
-\frac{b}{\lambda}(1-e^{-\lambda\left(  t-r\right)  })\text{ if }t\geq r,
\end{array}
\right.
\end{align*}
where $r\in\mathbb{R}$ is arbitrary. It is clear that $\mathcal{A}_{\infty
}=[-\frac{b}{\lambda},\frac{b}{\lambda}].$

\begin{lemma}
\label{AsympAutonOrdIncl}$\mathcal{R}$ is asymptotically autononomous.
\end{lemma}

\begin{proof}
Let $u_{n}\in\mathcal{R}_{\tau_{n}}$ be such that $\tau_{n}\rightarrow+\infty$
and $u\left(  \tau_{n}\right)  \rightarrow u_{0}.$

First, let $u_{0}\not =0$. For instance, take $u_{0}>0$. Hence, we can assume
that $u\left(  \tau_{n}\right)  >0$. In such a case, the solutions $u_{n}$ are
unique \cite[Corollary 3.5]{CLV16} and%
\[
u_{n}\left(  t\right)  =e^{-\lambda\left(  t-\tau_{n}\right)  }u_{n}\left(
\tau_{n}\right)  +\int_{\tau_{n}}^{t}e^{-\lambda\left(  t-s\right)  }b\left(
s\right)  ds.
\]
Also, the unique solution to the limit autononomous problem with $u\left(
0\right)  =u_{0}$ is given by%
\[
u\left(  t\right)  =e^{-\lambda t}u_{0}+\frac{b}{\lambda}\left(  1-e^{-\lambda
t}\right)  .
\]
Then for $v_{n}\left(  \text{\textperiodcentered}\right)  =u_{n}\left(
\text{\textperiodcentered}+\tau_{n}\right)  $ we have%
\begin{align*}
\left\vert v_{n}\left(  t\right)  -u\left(  t\right)  \right\vert  &
=\left\vert e^{-\lambda t}u_{n}\left(  \tau_{n}\right)  +\int_{\tau_{n}%
}^{t+\tau_{n}}e^{-\lambda\left(  t+\tau_{n}-s\right)  }b\left(  s\right)
ds-e^{-\lambda t}u_{0}-\frac{b}{\lambda}\left(  1-e^{-\lambda t}\right)
\right\vert \\
&  =\left\vert e^{-\lambda t}\left(  u_{n}\left(  \tau_{n}\right)
-u_{0}\right)  +\int_{\tau_{n}}^{t+\tau_{n}}e^{-\lambda\left(  t+\tau
_{n}-s\right)  }\left(  b\left(  s\right)  -b\right)  ds\right\vert \\
&  \leq e^{-\lambda t}\left\vert u_{n}\left(  \tau_{n}\right)  -u_{0}%
\right\vert +\frac{1}{\lambda}\sup_{s\geq\tau_{n}}\left\vert b\left(
s\right)  -b\right\vert \rightarrow0\text{ as }n\rightarrow\infty.
\end{align*}

Second, let $u_{0}=0$. If there is a subsequence $\{u\left(  \tau_{n_{k}%
}\right)  \}$ such that $u\left(  \tau_{n_{k}}\right)  >0$, then arguing as
before $v_{n_{k}}\left(  t\right)  $ converges to $u\left(  t\right)
=\frac{b}{\lambda}\left(  1-e^{-\lambda t}\right)  $, which is a solution to
the autonomous problem with $u\left(  0\right)  =0$. Assume then that
$u\left(  \tau_{n}\right)  =0$ for any $n$. Then the solutions $u_{n}$ are
either $u_{n}\left(  t\right)  \equiv0$ or of the form%
\begin{align*}
u_{n}^{+}\left(  t\right)   &  =\left\{
\begin{array}
[c]{c}%
0\text{ if }\tau_{n}\leq t\leq r_{n},\\
\int_{r_{n}}^{t}e^{-\lambda\left(  t-s\right)  }b\left(  s\right)  ds\text{ if
}t\geq r_{n},
\end{array}
\right. \\
u_{n}^{-}\left(  t\right)   &  =\left\{
\begin{array}
[c]{c}%
0\text{ if }\tau_{n}\leq t\leq r_{n},\\
-\int_{r_{n}}^{t}e^{-\lambda\left(  t-s\right)  }b\left(  s\right)  ds\text{
if }t\geq r_{n}.
\end{array}
\right.
\end{align*}
The case when $u_{n}\left(  t\right)  \equiv0$ (at least for a subsequence) is
trivial. Up to a subsequence let, for instance, $u_{n}=u_{n}^{+}$. Hence,%
\[
v_{n}(t)=\left\{
\begin{array}
[c]{c}%
0\text{ if }0\leq t\leq r_{n}-\tau_{n},\\
\int_{r_{n}}^{t+\tau_{n}}e^{-\lambda\left(  t+\tau_{n}-s\right)  }b\left(
s\right)  ds\text{ if }t\geq r_{n}-\tau_{n}.
\end{array}
\right.
\]
If $r_{n}-\tau_{n}\rightarrow+\infty$, then $v_{n}\left(  t\right)
\rightarrow0$ for all $t\geq0$. If not, then up to a subsequence $r_{n}%
-\tau_{n}\rightarrow\alpha_{0}$. Thus,%
\[
v_{n}(t)\rightarrow0\text{ if }0\leq t\leq\alpha_{0},
\]%
\[
v_{n}(t)=\int_{r_{n}}^{t+\tau_{n}}e^{-\lambda\left(  t+\tau_{n}-s\right)
}b\left(  s\right)  ds\rightarrow\frac{b}{\lambda}\left(  1-e^{-\lambda\left(
t-\alpha_{0}\right)  }\right)  \text{ if }\alpha_{0}\leq t.
\]
Therefore, $v_{n}\left(  t\right)  $ converges to a solution $u\left(
t\right)  $ to the autonomous problem with $u\left(  0\right)  =0.$
\end{proof}

\bigskip

\begin{lemma}
\label{ConvergFixedODE}$\lim_{t\rightarrow+\infty}\xi_{M}^{\pm}\left(
t\right)  =z_{1}^{\pm}.$
\end{lemma}

\begin{proof}
For $z_{1}^{+}$ we have%
\[
\xi_{M}^{+}\left(  t\right)  -z_{1}^{+}=\int_{-\infty}^{t}e^{-\lambda\left(
t-s\right)  }b\left(  s\right)  ds-\frac{b}{\lambda}=\int_{-\infty}%
^{t}e^{-\lambda\left(  t-s\right)  }\left(  b\left(  s\right)  -b\right)  ds.
\]
For any $\varepsilon>0$ we choose $t_{0}\left(  \epsilon\right)  $ such that
$\left\vert b\left(  s\right)  -b\right\vert \leq\frac{\varepsilon\lambda}{2}$
if $s\geq t_{0}$. Then we take $t_{1}\left(  \varepsilon\right)  $ such that
$2b_{1}e^{-\lambda t}e^{\lambda t_{0}}\leq\frac{\varepsilon\lambda}{2}$ for
$t\geq t_{1}$. Hence,%
\begin{align*}
\left\vert \xi_{M}^{+}\left(  t\right)  -z_{1}^{+}\right\vert  &  \leq
\int_{-\infty}^{t_{0}}e^{-\lambda\left(  t-s\right)  }\left\vert b\left(
s\right)  -b\right\vert ds+\int_{t_{0}}^{t}e^{-\lambda\left(  t-s\right)
}\left\vert b\left(  s\right)  -b\right\vert ds\\
&  \leq\frac{2b_{1}}{\lambda}e^{-\lambda t}e^{\lambda t_{0}}+\frac{1}{\lambda
}\sup_{s\geq t_{0}}\left\vert b\left(  s\right)  -b\right\vert \leq
\varepsilon,
\end{align*}
if $t\geq t_{1}.$ For $z_{1}^{-}$ the proof is analogous.
\end{proof}

\bigskip

\begin{lemma}
\label{CondAInfMOrdIncl}$\lim_{t\rightarrow+\infty}dist_{H}\left(
\mathcal{A}\left(  t\right)  ,\mathcal{A}_{\infty}\right)  =0.$
\end{lemma}

\begin{proof}
Let us prove that $\lim_{t\rightarrow+\infty}dist\left(  \mathcal{A}_{\infty
},\mathcal{A}\left(  t\right)  \right)  =0.$ If $z_{0}=0$, then $dist\left(
z_{0},\mathcal{A}\left(  t\right)  \right)  =0$ for any $t$. We take $z_{0}%
\in(0,\frac{b}{\lambda}]$ and $\varepsilon>0$ arbitrary. If $z_{0}\geq\frac
{b}{\lambda}-\frac{\varepsilon}{2}$, then Lemma \ref{ConvergFixedODE} implies
the existence of $T_{0}(\varepsilon)$ such that
\[
dist(z_{0},\mathcal{A}\left(  t\right)  )\leq\rho(z_{0},z_{1}^{+})+\rho\left(
z_{1}^{+},\xi_{M}^{+}\left(  t\right)  \right)  \leq\varepsilon\text{ for any
}t\geq T_{0}.
\]
If $z_{0}\leq\frac{b}{\lambda}-\frac{\varepsilon}{2}$, there is a unique
$a_{0}\in\lbrack0,a_{1}]$, where $a_{1}=\frac{1}{\lambda}\log\left(  \frac
{2b}{\lambda\varepsilon}\right)  ,$ such that
\[
z_{0}=\frac{b}{\lambda}\left(  1-e^{-a_{0}\lambda}\right)  .
\]
We take $T_{1}\left(  \epsilon\right)  $ such that $\left\vert b\left(
s\right)  -b\right\vert \leq\varepsilon\lambda$ if $s\geq T_{1}$. Then for any
$t\geq$ $T_{1}+a_{1}$ we obtain%
\begin{align*}
\left\vert z_{0}-u_{t-a_{0}}^{+}(t)\right\vert  &  =\left\vert \frac
{b}{\lambda}\left(  1-e^{-a_{0}\lambda}\right)  -\int_{t-a_{0}}^{t}%
e^{-\lambda\left(  t-s\right)  }b\left(  s\right)  ds\right\vert =\left\vert
\int_{t-a_{0}}^{t}e^{-\lambda\left(  t-s\right)  }(b-b\left(  s\right)
)ds\right\vert \\
&  \leq\frac{1}{\lambda}\sup_{s\geq t-a_{0}}\left\vert b-b\left(  s\right)
\right\vert \leq\frac{1}{\lambda}\sup_{s\geq t-a_{1}}\left\vert b-b\left(
s\right)  \right\vert \leq\varepsilon.
\end{align*}
For $z_{0}\in\lbrack-\frac{b}{\lambda},0)$ the proof is similar. Thus, we have
proved that for any $\varepsilon>0$ there is $T\left(  \varepsilon\right)  $
such that%
\[
dist\left(  \mathcal{A}_{\infty},\mathcal{A}\left(  t\right)  \right)
\leq\varepsilon\text{ if }t\geq T.
\]

The fact that $\lim_{t\rightarrow+\infty}dist\left(  \mathcal{A}\left(
t\right)  ,\mathcal{A}_{\infty}\right)  =0$ follows from Lemmas
\ref{ForwardACOrdIncl}, \ref{AsympAutonOrdIncl} and \ref{ConvergAutonM}.
\end{proof}

\bigskip

\begin{theorem}
$\mathcal{A}$ is a forward attractor.
\end{theorem}

\begin{proof}
It follows from Lemmas \ref{ForwardACOrdIncl}, \ref{AsympAutonOrdIncl},
\ref{CondAInfMOrdIncl} and Theorem \ref{ForwardAttractorTh2M}.
\end{proof}

\bigskip

\subsection{A parabolic differential inclusion}

We will study the problem
\begin{equation}
\left\{
\begin{array}
[c]{l}%
\dfrac{\partial u}{\partial t}-\dfrac{\partial^{2}u}{\partial x^{2}}\in
b(t)H_{0}(u)+\omega(t)u,\text{ on }(\tau,\infty)\times(0,1),\\
u(t,0)=u(t,1)=0,\\
u(\tau,x)=u_{\tau}(x),\text{ }x\in\left(  0,1\right)  ,
\end{array}
\right.  \label{Incl}%
\end{equation}
where $b:\mathbb{R}\rightarrow\mathbb{R}^{+},$ $\omega:\mathbb{R}%
\rightarrow\mathbb{R}^{+}$ are continuous functions such that
\begin{equation}
0<b_{0}\leq b\left(  t\right)  \leq b_{1},\ 0\leq\omega_{0}\leq\omega\left(
t\right)  \leq\omega_{1}, \label{Cond1}%
\end{equation}
and $H_{0}$ is the Heaviside function given in (\ref{Heaviside}).

Let $A:D(A)\rightarrow H,\ D(A)=H^{2}(0,1)\cap V,$ be the operator
$A=-\dfrac{d^{2}}{dx^{2}}$ with Dirichlet boundary conditions. This operator
is the generator of a $C_{0}$-semigroup $T(t)=e^{-At}$.

For any $u_{\tau}\in H$ the function $u\in C([\tau,+\infty),H)$ is said to be
a strong solution to problem (\ref{Incl}) if $u(\tau)=u_{\tau},\ u\left(
\text{\textperiodcentered}\right)  $ is absolutely continuous on $[T_{1}%
,T_{2}]$ for any $\tau<T_{1}<T_{2},$ $u\left(  t\right)  \in D(A)$ for a.a.
$t\in\left(  T_{1},T_{2}\right)  $, and there exists a function $r\in
L_{loc}^{2}([\tau,+\infty);H)$ such that $r\left(  t,x\right)  \in
b(t)H_{0}(u\left(  t,x\right)  )$ for a.a. $\left(  t,x\right)  \in\left(
\tau,+\infty\right)  \times\left(  0,1\right)  $ and%
\begin{equation}
\frac{du}{dt}+Au(t)=r(t)+\omega(t)u\text{ for a.a. }t\in\left(  \tau
,+\infty\right)  , \label{Equality}%
\end{equation}
where the equality is understood in the sense of the space $H.$

We will focus on non-negative solutions. Under assumption (\ref{Cond1}) it is
known \cite[Corollary 5]{CVL20} that for any $u_{\tau}\in H$ such that
$u_{\tau}\geq0$ there is at least one strong solution $u\left(
\text{\textperiodcentered}\right)  $ to problem (\ref{Incl}) such that
$u\left(  t\right)  \geq0$ for any $t\geq\tau$. The following facts were
proved in \cite{Valero2021}:

\begin{itemize}
\item Any strong solution $u$ to problem (\ref{Incl}) satisfies $u\in
C((\tau,+\infty),V).$

\item Let $u_{\tau}\in H$ be such that $u_{\tau}\geq0$ but $u_{\tau
}\not \equiv 0$ and let $u\left(  \text{\textperiodcentered}\right)  $ be a
non-negative solution to problem (\ref{Incl}). Then the solution $u\left(
\text{\textperiodcentered}\right)  $ is unique in the class of non-negative
solutions and $u(t)$ is positive for any $t>\tau.$ In addition, it is the
unique solution to the problem%
\begin{equation}
\left\{
\begin{array}
[c]{l}%
\dfrac{\partial u}{\partial t}-\dfrac{\partial^{2}u}{\partial x^{2}%
}=b(t)+\omega(t)u,\text{ on }(\tau,\infty)\times(0,1),\\
u(t,0)=u(t,1)=0,\\
u(\tau,x)=u_{\tau}(x),\text{ }x\in\left(  0,1\right)  .
\end{array}
\right.  \label{Aux2}%
\end{equation}

\item If $u_{\tau}\equiv0$, then, apart from the zero solution, all the
possible non-negative solutions are of the following type:%
\begin{equation}
u(t)=\left\{
\begin{array}
[c]{c}%
0\text{ if }\tau\leq t\leq t_{0},\\
u_{t_{0}}(t)\text{ if }t\geq t_{0},
\end{array}
\right.  \label{Sol0}%
\end{equation}
where $u_{t_{0}}($\textperiodcentered$)$ is the unique solution to the problem%
\begin{equation}
\left\{
\begin{array}
[c]{l}%
\dfrac{\partial u}{\partial t}-\dfrac{\partial^{2}u}{\partial x^{2}%
}=b(t)+\omega(t)u(t),\text{ on }(t_{0},\infty)\times(0,1),\\
u(t,0)=u(t,1)=0,\\
u(t_{0},x)=0,
\end{array}
\right.  \label{Aux3}%
\end{equation}
and $u(t)$ is positive for all $t>t_{0}.$

\item If $u_{\tau}\in V$ is such that $u_{\tau}\geq0$, $u_{\tau}\left(
x_{0}\right)  =0$ at some $x_{0}\in\left(  0,1\right)  $ but $u_{\tau
}\not \equiv 0$, then there cannot exist a non-negative solution backwards in time.

\item If $u_{\tau}\in H$ is such that $u_{\tau}\geq0$ and $u_{\tau}\not \in
V$, then there cannot exist a non-negative solution backwards in time.

\item If $u_{\tau}\equiv0$, then the unique non-negative solution backwards in
time is the zero solution, that is, $u\left(  t\right)  \equiv0$ for
$t\leq\tau$.
\end{itemize}

In order to study the pullback attractor we need to assume additionally that
\begin{equation}
\omega_{1}<\pi^{2}. \label{CondAttr2}%
\end{equation}

Let $H^{+}$ be the positive cone of $H$, that is,%
\[
H^{+}=\{v\in H:v(x)\geq0\text{ for a.a. }x\in\left(  0,1\right)  \}.
\]
We denote by $\mathcal{D}_{\tau}^{+}(u_{\tau})$ the set of all non-negative
solutions to problem (\ref{Incl}) with initial condition $u_{\tau}\in H^{+}$
at time $\tau$ and let $\mathcal{R}_{\tau}^{+}=\cup_{u_{\tau}\in H}%
\mathcal{D}_{\tau}^{+}(u_{\tau})$, $\mathcal{R}^{+}\mathcal{=\cup}_{\tau
\in\mathbb{R}}\mathcal{R}_{\tau}^{+}$. We define the map $U^{+}:\mathbb{R}%
_{\geq}\times H^{+}\rightarrow P(H^{+})$ given by%
\[
U^{+}(t,\tau,u_{\tau})=\{u(t):u\in\mathcal{D}_{\tau}^{+}(u_{\tau})\}.
\]
The family $\mathcal{R}^{+}$ satisfies the properties (K1)-(K4) \cite[Section
4]{CVL20}. We summarize several results from \cite{Valero2021}. The
corresponding strict multivalued process $U^{+}$ possesses a strictly
invariant pullback attractor $\mathcal{A}^{+}=\{\mathcal{A}^{+}(t)\}$.
Moreover,%
\[
\mathcal{A}^{+}(t)=\{\gamma(t):\gamma\text{ is a bounded complete trajectory
of }\mathcal{R}^{+}\},
\]
$\cup_{t\in\mathbb{R}}\mathcal{A}^{+}(t)$ is bounded in $V$, $\overline
{\cup_{t\in\mathbb{R}}\mathcal{A}^{+}(t)}$ is compact in $H$ and the sets
$\mathcal{A}^{+}(t)$ are compact in $V$.

The structure of the pullback attractor $\mathcal{A}^{+}$ is as follows.
First, there exists a bounded complete trajectory $\xi_{M}\left(  t\right)  $
such that:

\begin{enumerate}
\item $\xi_{M}\left(  t\right)  >0$ for all $t\in\mathbb{R}$. Hence, for any
$\tau\in\mathbb{R}$ it is the unique solution to (\ref{Aux2}) with $u_{\tau
}=\xi_{M}(\tau).$

\item $0\leq\gamma\left(  t\right)  \leq\xi_{M}\left(  t\right)  $ for any
bounded complete trajectory $\gamma$ of $\mathcal{R}^{+};$

\item $\xi_{M}$ is the unique bounded complete trajectory of $\mathcal{R}^{+}$
such that $\xi_{M}\left(  t\right)  >0$ for all $t\in\mathbb{R};$
\end{enumerate}

The solution $\xi_{M}$ is a so called nonautonomous equilibrium. Second, any
bounded complete trajectory $\gamma$ of $\mathcal{R}^{+}$ distinct from $0$
and $\xi_{M}$ has the form%
\begin{equation}
\gamma\left(  t\right)  =\left\{
\begin{array}
[c]{c}%
0\text{ if }t\leq t_{0},\\
u\left(  t\right)  \text{ if }t\geq t_{0},
\end{array}
\right.  \label{CompleteTrayIncl}%
\end{equation}
for some $t_{0}\in\mathbb{R}$, where $u$ is the unique solution to
(\ref{Aux3}). Third,%
\[
\left\Vert \gamma\left(  t\right)  -\xi_{M}\left(  t\right)  \right\Vert
\rightarrow0\text{ as }t\rightarrow+\infty.
\]
Thus, the pullback attractor consists of the equilibria $0$, $\xi_{M}$ and the
heteroclinic connections between them.

\begin{lemma}
\label{ACIncl}$U^{+}$ is forward asymptotically compact.
\end{lemma}

\begin{proof}
Let $y_{n}\in U^{+}(t_{n},\tau,B)$, where $t_{n}\rightarrow+\infty$ and $B$ is
a bounded set. Then $y_{n}=u^{n}\left(  t_{n}\right)  $, where $u^{n}$ are
strong solutions to (\ref{Incl}) with $u^{n}(\tau)=u_{\tau}^{n}\in B.$ Then
there exist $f^{n}\in L_{loc}^{2}(\tau,+\infty;H)$ such that $f^{n}(t,x)\in
H_{0}(u^{n}(t,x))$ a.e. in $(\tau,+\infty)\times\left(  0,1\right)  $ and%
\begin{equation}
\dfrac{du^{n}}{dt}(t)-Au^{n}(t)=b(t)f^{n}(t)+\omega(t)u^{n}(t),\quad\text{a.e.
in}\;(\tau,+\infty), \label{P0}%
\end{equation}
Multiplying by $u^{n}$ we have%
\[
\frac{1}{2}\frac{d}{dt}\left\Vert u^{n}\right\Vert ^{2}+\Vert u^{n}\Vert
_{V}^{2}\leq b_{1}\int_{0}^{1}\left\vert u^{n}\right\vert dx+\omega
_{1}\left\Vert u^{n}\right\Vert ^{2}\leq C_{1}+\left(  \omega_{1}+\delta
_{0}\right)  \Vert u^{n}\Vert^{2},
\]
where $\delta_{0}=\frac{\pi^{2}-\omega_{1}}{2}$, and $\Vert z\Vert_{V}^{2}%
\geq\pi^{2}\left\Vert z\right\Vert ^{2}$ implies%
\[
\left\Vert u^{n}\left(  t\right)  \right\Vert ^{2}\leq e^{-(\pi^{2}-\omega
_{1})(t-\tau)}\left\Vert u_{\tau}^{n}\right\Vert +\frac{2C_{1}}{\pi^{2}%
-\omega_{1}}\text{ for all }t\geq\tau.
\]
Then there exists $T_{1}>0$ such that%
\[
\left\Vert u^{n}\left(  t\right)  \right\Vert ^{2}\leq1+\frac{2C_{1}}{\pi
^{2}-\omega_{1}}=C_{2}\text{ for all }t\geq T_{1}+\tau.
\]
Therefore,
\begin{equation}
\int_{t}^{t+1}\Vert u^{n}\Vert_{V}^{2}ds\leq C_{1}+\frac{C_{2}}{2}+\left(
\omega_{1}+\delta_{0}\right)  C_{2}=C_{3}\text{ if }t\geq T_{1}+\tau.
\label{AcotV}%
\end{equation}
Multiplying (\ref{P0}) by $-\frac{\partial^{2}u}{\partial x^{2}}$ and using
Young's inequality we obtain
\[
\frac{d}{dt}\Vert u^{n}\Vert_{V}^{2}+2\Vert\frac{\partial^{2}u^{n}}{\partial
x^{2}}\left(  r\right)  \Vert^{2}\leq2b_{1}^{2}\Vert f^{n}\left(  r\right)
\Vert^{2}+2\omega_{1}^{2}\left\Vert u^{n}\left(  r\right)  \right\Vert
^{2}+\Vert\frac{\partial^{2}u^{n}}{\partial x^{2}}\left(  s\right)  \Vert
^{2}.
\]
These operations are correct (see, for instance, \cite[p.1070]{CCMV19}).
Integrating over $\left(  s,t\right)  $ with $T_{1}+\tau\leq t-1<s<t$,we have%
\[
\Vert u^{n}(t)\Vert_{V}^{2}\leq\Vert u^{n}(s)\Vert_{V}^{2}+2b_{1}^{2}\int%
_{s}^{t}\Vert f^{n}\left(  r\right)  \Vert^{2}dr+2\omega_{1}^{2}C_{2}.
\]
Integrating over $\left(  t-1,t-1+\varepsilon\right)  $ with $0<\varepsilon
<1$, and using (\ref{AcotV}) we infer%
\[
\Vert u^{n}(t)\Vert_{V}^{2}\leq\frac{C_{4}}{\varepsilon}\text{ if }%
t\geq1+T_{1}+\tau.
\]

Hence, the sequence $\{y_{n}\}$ is bounded in $V$. The compact embedding
$V\subset H$ implies that $\{y_{n}\}$ is relatively compact in $H.$
\end{proof}

\begin{remark}
This result is valid for the process $U$ generated by all the solutions, not
only the non-negative ones.
\end{remark}

\begin{theorem}
$\mathcal{A}^{+}$ is a forward attractor.
\end{theorem}

\begin{proof}
Let $B$ be a bounded set and $R>0$ such that $\left\vert z\right\vert \leq R$
for any $z\in B$.

If $u_{\tau}\not =0$, then $U^{+}(t,\tau,u_{\tau})=u\left(  t\right)  $, where
$u$ is the unique solution to (\ref{Aux2}). The difference $w=u-\xi_{M}$
satisfies%
\[
\frac{dw}{dt}-w_{xx}=\omega(t)w.
\]
Hence,%
\[
\left\Vert w\left(  t\right)  \right\Vert ^{2}\leq e^{-2\left(  \pi^{2}%
-\omega_{1}\right)  (t-\tau)}\left\Vert w\left(  \tau\right)  \right\Vert
^{2}\leq2e^{-2\left(  \pi^{2}-\omega_{1}\right)  (t-\tau)}\left(  R^{2}%
+D^{2}\right)  \rightarrow0\text{ as }t\rightarrow+\infty,
\]
where $D$ is a uniform bound of $\left\Vert \xi_{M}\left(  t\right)
\right\Vert $.

If $u_{\tau}=0$, then the possible solutions are either $u\left(  t\right)
\equiv0$ or (\ref{Sol0}). Hence, $U^{+}(t,\tau,0)\subset\mathcal{A}^{+}\left(
t\right)  $ for all $t\geq\tau$.

We deduce that%
\[
dist(U(t,\tau,B),\mathcal{A}^{+}\left(  t\right)  )\rightarrow0\text{ as
}t\rightarrow+\infty.
\]
Thus, $\mathcal{A}^{+}$ is a forward attractor.
\end{proof}

\bigskip

As in the previous application, in the particular situation where problem
(\ref{Incl}) is asymptotically autonomous, we will prove that
(\ref{CondOmega0M}) is satisfied. We assume that
\[
b\left(  t\right)  \rightarrow b\in\lbrack b_{0},b_{1}]\ \omega\left(
t\right)  \rightarrow\omega\in\lbrack\omega_{0},\omega_{1}]\text{ as
}t\rightarrow+\infty.
\]
For the autonomous system with $b\left(  t\right)  \equiv b$,\ $\omega\left(
t\right)  \equiv\omega$ it is known \cite{CVL20} that if $\mathcal{R}_{0}^{+}$
is the set of all nonnegative strong solutions to (\ref{Incl}), then
properties (H1)-(H4) are satisfied. We summarize several results from
\cite{Valero2021}. The corresponding strict multivalued semiflow $G$ has the
strict invariant attractor $\mathcal{A}_{\infty}^{+}$, which is compact in
$V$. Moreover,
\[
\mathcal{A}_{\infty}^{+}=\{\gamma(t):\gamma\text{ is a bounded complete
trajectory of }\mathcal{R}_{0}^{+}\}.
\]
There is a positive fixed point $v_{1}^{+}$ and $0\leq y\leq v_{1}^{+}$ for
all $y\in\mathcal{A}_{\infty}^{+}$. The only bounded complete trajectories
distinct form the fized points $0,v_{1}^{+}$ are of the type
(\ref{CompleteTrayIncl}). They satisfy that%
\begin{align*}
\gamma\left(  t\right)   &  \rightarrow v_{1}^{+}\text{ as }t\rightarrow
+\infty,\\
\gamma\left(  t\right)   &  \rightarrow0\text{ as }t\rightarrow-\infty.
\end{align*}
Hence, the attractor consists of the fixed points $0,v_{1}^{+}$ an the
heteroclinic connections between them of the type (\ref{CompleteTrayIncl}).

Let us prove first that $\mathcal{R}^{+}$ is asymptotically autonomous. We
give a more general result that is valid for all the solutions, not only the
non-negative ones.

\begin{theorem}
\label{AsymptAutonIncl}Let $\tau_{n}\nearrow+\infty.$ If $u_{\tau_{n}}$ and
$u_{\tau_{n}}\rightarrow u_{0}$ in $L^{2}(0,1)$ as $n$ $\rightarrow$
$+\infty,$ then for each family of strong solutions $u^{n}$ of problem
(\ref{Incl}) with $u^{n}\left(  \tau_{n}\right)  =u_{\tau_{n}}$ there exists a
strong solution $v$ of the autonomous problem such that, up to a subsequence,
$v^{n}(t):=u^{n}(t+\tau_{n})$ $\rightarrow$ $v(t)$ in $H,$ as $n$
$\rightarrow$ $+\infty,$ uniformly on compact sets of $[0,+\infty).$
\end{theorem}

\begin{proof}
Let $u^{n}$ be strong solutions of (\ref{Incl}) with $u^{n}(\tau_{n}%
)=u_{\tau_{n}}.$ Then there exist $f^{n}\in L_{loc}^{2}(\tau_{n},+\infty;H)$
such that $f^{n}(t,x)\in H_{0}(u^{n}(t,x))$ a.e. in $(\tau_{n},+\infty
)\times\left(  0,1\right)  $ and%
\begin{equation}
\dfrac{du^{n}}{dt}(t)-Au^{n}(t)=b(t)f^{n}(t)+\omega(t)u^{n}(t),\quad\text{a.e.
in}\;(\tau_{n},+\infty), \label{P}%
\end{equation}

Multiplying by $u^{n}$ we have%
\[
\frac{1}{2}\frac{d}{dt}\left\Vert u^{n}\right\Vert ^{2}+\Vert u^{n}\Vert
_{V}^{2}\leq b_{1}\int_{0}^{1}\left\vert u^{n}\right\vert dx+\omega
_{1}\left\Vert u^{n}\right\Vert ^{2}\leq C_{1}+\left(  \omega_{1}+\delta
_{0}\right)  \Vert u^{n}\Vert^{2},
\]
where $\delta_{0}=\frac{\pi^{2}-\omega_{1}}{2}$. Since $\Vert z\Vert_{V}%
^{2}\geq\pi^{2}\left\Vert z\right\Vert ^{2}$, we deduce that%
\begin{equation}
\left\Vert u^{n}\left(  \tau_{n}+t\right)  \right\Vert ^{2}\leq e^{-(\pi
^{2}-\omega_{1})t}+\frac{2C_{1}}{\pi^{2}-\omega_{1}}\leq C_{2}\text{ for all
}t\geq0. \label{Acot1}%
\end{equation}
Let us fix an arbitrary $T>0$. Thus,
\begin{equation}
\int_{\tau_{n}}^{\tau_{n}+T}\Vert u^{n}\Vert_{V}^{2}ds\leq C_{1}T+\frac{C_{2}%
}{2}+\left(  \omega_{1}+\delta_{0}\right)  TC_{2}=C_{3,T}. \label{Acot2}%
\end{equation}
Multiplying (\ref{P}) by $-\frac{\partial^{2}u}{\partial x^{2}}$ and using
Young's inequality we obtain
\begin{equation}
\frac{d}{dt}\Vert u^{n}\Vert_{V}^{2}+2\Vert\frac{\partial^{2}u^{n}}{\partial
x^{2}}\left(  r\right)  \Vert^{2}\leq2b_{1}^{2}\Vert f^{n}\left(  r\right)
\Vert^{2}+2\omega_{1}^{2}\left\Vert u^{n}\left(  r\right)  \right\Vert
^{2}+\Vert\frac{\partial^{2}u^{n}}{\partial x^{2}}\left(  s\right)  \Vert^{2}.
\label{Ineq1}%
\end{equation}
These operations are correct (see, for instance, \cite[p.1070]{CCMV19}).
Integrating over $\left(  s,t+\tau_{n}\right)  $ with $\tau_{n}<s<t+\tau_{n}$
we have%
\[
\Vert u^{n}(t+\tau_{n})\Vert_{V}^{2}\leq\Vert u^{n}(s)\Vert_{V}^{2}+2b_{1}%
^{2}\int_{\tau_{n}}^{\tau_{n}+T}\Vert f^{n}\left(  r\right)  \Vert
^{2}dr+2\omega_{1}^{2}C_{2}T.
\]
Integrating over $\left(  \tau_{n},\tau_{n}+\varepsilon\right)  $ with
$0<\varepsilon\leq t$ and using (\ref{Acot2}) we infer%
\begin{equation}
\Vert u^{n}(t+\tau_{n})\Vert_{V}^{2}\leq\frac{C_{4,T}}{\varepsilon}\text{ for
}t\geq\varepsilon. \label{Acot3}%
\end{equation}
Thus, from (\ref{P}), (\ref{Ineq1}) and (\ref{Acot3}) we find that%
\begin{equation}
\int_{\tau_{n}+\varepsilon}^{\tau_{n}+T}\Vert\frac{\partial^{2}u^{n}}{\partial
x^{2}}\left(  s\right)  \Vert^{2}ds\leq C_{5,T}\text{,} \label{Ineq2}%
\end{equation}%
\begin{equation}
\int_{\tau_{n}+\varepsilon}^{\tau_{n}+T}\Vert\frac{du^{n}}{ds}\left(
s\right)  \Vert^{2}ds\leq C_{6,T}. \label{Ineq3}%
\end{equation}

From the definition of $H_{0}$ it follows that $\left\vert f^{n}\left(
t,x\right)  \right\vert \leq1$, so in particular the sequence $\{g^{n}\}$,
defined by $g^{n}($\textperiodcentered$)=f^{n}($\textperiodcentered$+\tau
_{n}),$ is bounded in $L^{\infty}(0,T;L^{2}(0,1))$. Hence, up to a
subsequence, $g^{n}\rightarrow g$ weakly star in $L^{\infty}(0,T;L^{2}(0,1))$
and weakly in $L^{2}(0,T;L^{2}(0,1))$ for some function $g.$ Let $v^{n}\left(
\text{\textperiodcentered}\right)  =u^{n}\left(  \text{\textperiodcentered
}+\tau_{n}\right)  $, $b^{n}\left(  \text{\textperiodcentered}\right)
=b\left(  \text{\textperiodcentered}+\tau_{n}\right)  ,\ \omega^{n}\left(
\text{\textperiodcentered}\right)  =\omega\left(  \text{\textperiodcentered
}+\tau_{n}\right)  $. Then for each $n$ the function $v^{n}$ is the unique
strong solution of the problem%
\[
\left\{
\begin{array}
[c]{l}%
\dfrac{\partial v^{n}}{\partial t}-\dfrac{\partial^{2}v^{n}}{\partial x^{2}%
}=b^{n}(t)g^{n}(t)+\omega^{n}(t)v^{n}(t),\text{ on }(0,T)\times\left(
0,1\right)  ,\\
v^{n}(t,0)=v^{n}(t,1)=0,\\
v^{n}(0,x)=u_{\tau_{n}}(x).
\end{array}
\right.
\]
Then (\ref{Acot1}), (\ref{Acot2}), (\ref{Acot3}), (\ref{Ineq2}), (\ref{Ineq3})
imply the existence of a function $v$ and a subsequence of $\{v^{n}\}$ such
that%
\[
v^{n}\rightarrow v\text{ weakly star in }L^{\infty}(0,T;H),
\]%
\[
v^{n}\rightarrow v\text{ weakly in }L^{2}(0,T;V),\text{ }%
\]%
\[
v^{n}\rightarrow v\text{ weakly star in }L^{\infty}(\varepsilon,T;V),
\]%
\[
v^{n}\rightarrow v\text{ weakly in }L^{2}(\varepsilon,T;D(A)),\text{ }%
\]%
\[
\frac{dv^{n}}{dt}\rightarrow\frac{dv}{dt}\text{ weakly in }L^{2}%
(\varepsilon,T;H),
\]
for all $\varepsilon>0$. The functions $v^{n}:[\varepsilon,T]\rightarrow H$
are then equicontinuous. As $v^{n}\left(  t\right)  $ is relatively compact in
$H$ for each $t\in\lbrack0,T]$, the Ascoli-Arzel\`{a} theorem gives that%
\[
v^{n}\rightarrow v\text{ in }C([\varepsilon,T],L^{2}(0,1)).
\]
Also, the convergences $b\left(  t\right)  \rightarrow b,\ \omega\left(
t\right)  \rightarrow\omega,$ as $t\rightarrow+\infty,$ imply that
\[
b^{n}\rightarrow b,\ \omega^{n}\rightarrow\omega\text{ in }C([0,T]).
\]
Passing to the limit in (\ref{P}) we obtain that%
\[
\frac{dv}{dt}-Av=bg+\omega v\text{ in }L^{2}(\varepsilon,T;H),
\]
so%
\[
\dfrac{dv}{dt}(t)-Av(t)=bg(t)+\omega v(t)\text{ in }H\text{ for a.a. }%
t\in\left(  0,T\right)  .
\]
Also, $v:[\varepsilon,T]\rightarrow L^{2}(0,1)$ is absolutely continuous.

We will verify that $g\left(  t,x\right)  \in H_{0}(v(t,x))$ for a.a. $\left(
t,x\right)  .$ For a.a. $\left(  t,x\right)  \in\left(  0,T\right)
\times\left(  0,1\right)  $ there exists $N(t,x)$ such that $g^{n}(t,x)\in
H_{0}(v(t,x))$ for any $n\geq N.$ Indeed, let $A\subset\lbrack0,T]\times
\lbrack0,1]$ be a set of zero measure such that $v^{n}(t,x)\rightarrow v(t,x)$
for any $\left(  t,x\right)  \in A^{c}$. If $(t_{0},x_{0})\in A^{c}$ satisfies
$v\left(  t_{0},x_{0}\right)  =0$, then the result follows from $g^{n}%
(t_{0},x_{0})\in\lbrack-1,1]=H_{0}(v(t_{0},x_{0}))$ for all $n.$ If
$(t_{0},x_{0})\in A^{c}$ is such that $v(t_{0},x_{0})>0$, then $v^{n}%
(t_{0},x_{0})\rightarrow v(t_{0},x_{0})$ implies the existence of
$N(t_{0},x_{0})$ for which $v^{n}(t_{0},x_{0})>0$ for $n\geq N$, and
consequently $g^{n}(t_{0},x_{0})=1=H_{0}(v(t_{0},x_{0}))$. The same argument
is valid for $v(t_{0},x_{0})<0$. As $g^{n}\rightarrow g$ weakly in
$L^{1}(0,T;H)$ and the set $H_{0}(v(t,x))$ is convex, we obtain from
\cite[Lemma 32]{Valero2021} that $g\left(  t,x\right)  \in H_{0}(v(t,x))$ for
a.a. $\left(  t,x\right)  $.

In order to show that $v$ is a strong solution it remains to prove that
$v\left(  t\right)  $ is continuous as $t\rightarrow0^{+}$. Let $z\left(
t\right)  $ be the unique solution to the problem
\[
\left\{
\begin{array}
[c]{l}%
\dfrac{\partial z}{\partial t}-\dfrac{\partial^{2}z}{\partial x^{2}}=\omega
z,\text{ on }(0,T)\times\left(  0,1\right)  ,\\
z(t,0)=z(t,1)=0,\\
z(0,x)=u_{0}(x).
\end{array}
\right.
\]
The difference $w_{n}\left(  t\right)  =v^{n}(t)-z\left(  t\right)  $
satisfies%
\begin{align*}
\frac{dw_{n}}{dt}-\dfrac{\partial^{2}w_{n}}{\partial x^{2}}  &  =b^{n}%
(t)g^{n}(t)+\omega^{n}\left(  t\right)  v^{n}-\omega z\\
&  =b^{n}(t)g^{n}(t)+\omega^{n}\left(  t\right)  w_{n}+\left(  \omega
^{n}\left(  t\right)  -\omega\right)  z.
\end{align*}
Using $b^{n}(t)\left\vert g^{n}(t,x)\right\vert \leq b_{1},\ \omega_{n}\left(
t\right)  \leq\omega_{1}$ we have%
\[
\frac{1}{2}\frac{d}{dt}\left\Vert w_{n}\right\Vert ^{2}+\frac{\pi^{2}%
-\omega_{1}}{2}\left\Vert w_{n}\right\Vert ^{2}\leq C(1+\left\vert \omega
^{n}\left(  t\right)  -\omega\right\vert ^{2}\left\Vert z\right\Vert ^{2}).
\]
We take $R_{0}>0$ such that $\left\Vert z\left(  t\right)  \right\Vert \leq
R_{0}$ for $t\in\lbrack0,T]$. Then%
\[
\left\Vert w_{n}(t)\right\Vert ^{2}\leq\left\Vert w_{n}\left(  0\right)
\right\Vert ^{2}+Ct+R_{0}^{2}t\sup_{s\in\lbrack0,t]}\left\vert \omega
^{n}\left(  s\right)  -\omega\right\vert ^{2}.
\]
Passing to the limit as $n\rightarrow\infty$ we find that%
\[
\left\Vert v\left(  t\right)  -z\left(  t\right)  \right\Vert ^{2}\leq
Ct\text{ for }t>0.
\]
Hence,
\[
\left\Vert v\left(  t\right)  -u_{0}\right\Vert \leq\left\Vert v\left(
t\right)  -z\left(  t\right)  \right\Vert +\left\Vert z\left(  t\right)
-u_{0}\right\Vert \leq\sqrt{Ct}+\left\Vert z\left(  t\right)  -u_{0}%
\right\Vert \rightarrow0\text{ as }t\rightarrow0^{+}.
\]

Finally, we need to check that $v^{n}\left(  t_{n}\right)  \rightarrow u_{0}$
as $t_{n}\rightarrow0$. This follows from%
\begin{align*}
\left\Vert v^{n}\left(  t_{n}\right)  -u_{0}\right\Vert  &  \leq\left\Vert
w\left(  t_{n}\right)  \right\Vert +\left\Vert z\left(  t_{n}\right)
-u_{0}\right\Vert \\
&  \leq\sqrt{\left\Vert w_{n}\left(  0\right)  \right\Vert ^{2}+Ct_{n}%
+R_{0}^{2}t_{n}\sup_{s\in\lbrack0,t_{n}]}\left\vert \omega^{n}\left(
s\right)  -\omega\right\vert ^{2}}+\left\Vert z\left(  t_{n}\right)
-u_{0}\right\Vert \rightarrow0\text{ as }n\rightarrow\infty.
\end{align*}
Thus, $v^{n}\rightarrow v$ in $C([0,T],H).$

By a diagonal argument we obtain the desired strong solution $v$ defined in
$[0,+\infty)$ and that%
\[
u^{n}(t+\tau_{n})=v^{n}(t)\rightarrow v(t)
\]
uniformly on compact sets of $[0,+\infty).$
\end{proof}

\bigskip

The following facts were established in \cite{LSSV};

\begin{itemize}
\item $\xi_{M}\left(  t\right)  \rightarrow v_{1}^{+}$ as $t\rightarrow
+\infty.$

\item $\lim_{t\rightarrow+\infty}dist\left(  \mathcal{A}^{+}(t),\mathcal{A}%
_{\infty}^{+}\right)  =0.$
\end{itemize}

\begin{remark}
In \cite{LSSV} it was proved that $\lim_{t\rightarrow+\infty}dist\left(
\mathcal{A}(t),\mathcal{A}_{\infty}\right)  =0$, where the attractors refer to
the whole set of solutions, not only the non-negative ones. The convergence
$\lim_{t\rightarrow+\infty}dist\left(  \mathcal{A}^{+}(t),\mathcal{A}_{\infty
}^{+}\right)  =0$ follows directly from this more general result. This result
follows also from Lemma \ref{ACIncl}, Theorem \ref{AsymptAutonIncl} and Lemma
\ref{ConvergAutonM}.
\end{remark}

Let us establish condition (\ref{CondAInfM}).

\begin{lemma}
\label{CondAInfMIncl}$\lim_{t\rightarrow+\infty}dist_{H}\left(  \mathcal{A}%
^{+}\left(  t\right)  ,\mathcal{A}_{\infty}^{+}\right)  =0.$
\end{lemma}

\begin{proof}
We only need to prove that $\lim_{t\rightarrow+\infty}dist\left(
\mathcal{A}_{\infty}^{+},\mathcal{A}^{+}\left(  t\right)  \right)  =0.$

If $z_{0}=0$, then $dist\left(  z_{0},\mathcal{A}^{+}\left(  t\right)
\right)  =0$ for all $t$.

Let $\varepsilon>0$ be arbitrary. If $z_{0}\in\mathcal{A}_{\infty}^{+}$ is
such that $\rho\left(  z_{0},v_{1}^{+}\right)  \leq\frac{\varepsilon}{2}$,
then $\xi_{M}\left(  t\right)  \rightarrow v_{1}^{+}$ implies that there is
$T_{0}\left(  \varepsilon\right)  $ for which
\[
dist(z_{0},\mathcal{A}^{+}\left(  t\right)  )\leq\left\Vert z_{0}-z_{1}%
^{+}\right\Vert +\left\Vert z_{1}^{+}-\xi_{M}^{+}\left(  t\right)  \right\Vert
\leq\varepsilon\text{ for any }t\geq T_{0}.
\]
Let now $\rho\left(  z_{0},v_{1}^{+}\right)  >\frac{\varepsilon}{2}$,
$z_{0}\not =0$. There exists $a_{0}>0$ such that $z_{0}=\varphi_{0}\left(
a_{0}+t_{0}\right)  $, where $\varphi_{0}$ is a bounded complete trajectry of
the type (\ref{CompleteTrayIncl}) with $b\left(  t\right)  =b,\ \omega\left(
t\right)  =\omega$ and $t_{0}$ is arbitrary. Since $\varphi_{0}\left(
t\right)  \overset{t\rightarrow+\infty}{\rightarrow}v_{1}^{+}$, there is
$a_{1}\left(  \varepsilon\right)  >0$ such that $\left\Vert \varphi_{0}\left(
t\right)  -z_{1}^{+}\right\Vert \leq\frac{\varepsilon}{2}$ if $t\geq
a_{1}+t_{0}$. Hence, $a_{0}<a_{1}$. We pick then the bounded complete
trajectory $\varphi_{t-a_{0}}$ of the type (\ref{CompleteTrayIncl}) for the
nonautonomous problem with $t_{0}=t-a_{0}$. Then the difference $w=\varphi
_{t-a_{0}}-\varphi_{0}$ satisfies%
\begin{align*}
\frac{dw}{ds}-w_{xx}  &  =b\left(  s\right)  -b+\omega\left(  s\right)
\varphi_{t-a_{0}}-\omega\varphi_{0}\\
&  =b\left(  s\right)  -b+\omega\left(  s\right)  w+\left(  \omega\left(
s\right)  -\omega\right)  \varphi_{0}.
\end{align*}
As $\varphi_{0}\left(  s\right)  \in\mathcal{A}_{\infty}^{+}$, there is $R>0$
such that $\left\Vert \varphi_{0}\left(  s\right)  \right\Vert \leq R$ for all
$s\in\mathbb{R}$. Hence,%
\[
\frac{1}{2}\frac{d}{dt}\left\Vert w\right\Vert ^{2}+\frac{\pi^{2}-\omega_{1}%
}{2}\left\Vert w\left(  t\right)  \right\Vert ^{2}\leq C(\left\vert b\left(
s\right)  -b\right\vert ^{2}+\left\vert \omega\left(  s\right)  -\omega
\right\vert ^{2}).
\]
Thus,%
\[
\left\Vert w\left(  t\right)  \right\Vert ^{2}\leq e^{-\left(  \pi^{2}%
-\omega_{1}\right)  a_{0}}\left\Vert w\left(  t_{0}\right)  \right\Vert
^{2}+2C\int_{t-a_{0}}^{t}e^{-\left(  \pi^{2}-\omega_{1}\right)  (t-s)}\left(
\left\vert b\left(  s\right)  -b\right\vert ^{2}+\left\vert \omega\left(
s\right)  -\omega\right\vert ^{2}\right)  ds.
\]
We choose $T_{1}\left(  \varepsilon\right)  >0$ such that%
\[
\left\vert b\left(  s\right)  -b\right\vert ^{2}+\left\vert \omega\left(
s\right)  -\omega\right\vert ^{2}\leq\frac{\varepsilon^{2}\left(  \pi
^{2}-\omega_{1}\right)  }{2C}\text{ for any }s\geq T_{1}.
\]
Then $w\left(  t_{0}\right)  =0$ and $w\left(  t\right)  =\varphi_{t-a_{0}%
}\left(  t\right)  -z_{0}$ imply that if $t\geq T_{1}+a_{1}$ we obtain%
\begin{align*}
\left\Vert \varphi_{t-a_{0}}\left(  t\right)  -z_{0}\right\Vert ^{2}  &
\leq\frac{2C}{\pi^{2}-\omega_{1}}\sup_{s\geq t-a_{0}}\left(  \left\vert
b\left(  s\right)  -b\right\vert ^{2}+\left\vert \omega\left(  s\right)
-\omega\right\vert ^{2}\right) \\
&  \leq\frac{2C}{\pi^{2}-\omega_{1}}\sup_{s\geq t-a_{1}}\left(  \left\vert
b\left(  s\right)  -b\right\vert ^{2}+\left\vert \omega\left(  s\right)
-\omega\right\vert ^{2}\right)  \leq\varepsilon.
\end{align*}

Therefore, we have proved that for any $\varepsilon>0$ there is $T\left(
\varepsilon\right)  $ such that%
\[
dist\left(  \mathcal{A}_{\infty}^{+},\mathcal{A}^{+}\left(  t\right)  \right)
\leq\varepsilon\text{ if }t\geq T.
\]

\end{proof}

\bigskip

\begin{theorem}
$\mathcal{A}^{+}$ is a forward attractor.
\end{theorem}

\begin{proof}
It follows from Lemmas \ref{ACIncl}, \ref{CondAInfMIncl} and Theorems
\ref{AsymptAutonIncl}, \ref{ForwardAttractorTh2M}.
\end{proof}

\bigskip

\textbf{Acknowledgments}

The author has been supported by Ministerio de Ciencia, Innovaci\'{o}n y
Universidades, AEI and FEDER, grant PID2024-156228NB-I00, and by Generalitat
Valenciana, Conselleria de Educaci\'{o}n, Cultura, Universidades y Empleo,
grant CIAICO/2024/251\textbf{.}

\end{document}